\documentstyle[11pt]{article}
       \font\tenmsb=msbm10
       \font\sevenmsb=msbm7
       \font\fivemsb=msbm5
       \catcode`\@=11
       \ifx\amstexloaded@\relax\catcode`\@=\active
       \endinput\else\let\amstexloaded@\relax\fi
       \def\spaces@{\space\space\space\space\space}
       \def\spaces@@{\spaces@\spaces@\spaces@\spaces@\spaces@}
       \def\space@.{\futurelet\space@\relax}
       \space@. %
       \def\Err@#1{\errhelp\defaulthelp@\errmessage{AmS-TeX error: #1}}
       \def\relaxnext@{\let\next\relax}
       \def\accentfam@{7}
       \def\noaccents@{\def\accentfam@{0}}
       \def\Cal{\relaxnext@\ifmmode\let\next\Cal@\else
       \def\next{\Err@{Use \string\Cal\space only in math mode}}\fi\next}
       \def\Cal@#1{{\Cal@@{#1}}}
       \def\Cal@@#1{\noaccents@\fam\tw@#1}
       \def\Bbb{\relaxnext@\ifmmode\let\next\Bbb@\else
       \def\next{\Err@{Use \string\Bbb\space only in math mode}}\fi\next}
       \def\Bbb@#1{{\Bbb@@{#1}}}
       \def\Bbb@@#1{\noaccents@\fam\msbfam#1}
       \newfam\msbfam
       \textfont\msbfam=\tenmsb
       \scriptfont\msbfam=\sevenmsb
       \scriptscriptfont\msbfam=\fivemsb

\def\notin{ \in \! \!\!\!\!\!  / }

\def\Z{{\Bbb Z}}
\def\R{{\Bbb R}}
\def\C{{\Bbb C}}
\def\Q{{\Bbb Q}}

\newtheorem{Theorem}{Theorem}
\newtheorem{Lemma}{Lemma}[section]
\newtheorem{Proposition}{Proposition}[section]
\newtheorem{Corollary}{Corollary}[section]

\newtheorem{Definition}{Definition}[section]

\@addtoreset{equation}{}

\newcommand{\qed}{\nolinebreak\hfill$\Box$
\par\medbreak}

\newcommand{\bq}{\begin{equation}}
\newcommand{\eq}{\end{equation}}
\newcommand{\bp}{\begin{Proposition}}
\newcommand{\ep}{\end{Proposition}}
\newcommand{\bdf}{\begin{Definition}}
\newcommand{\edf}{\end{Definition}}
\newcommand{\bl}{\begin{Lemma}}
\newcommand{\el}{\end{Lemma}}
\newcommand{\ba}{\begin{array}}
\newcommand{\ea}{\end{array}}
\newcommand{\bea}{\begin{eqnarray}}
\newcommand{\eea}{\end{eqnarray}}

\begin{document}
\setlength{\columnsep}{5pt}
\title{\bf Floer Homology for Symplectomorphism}
\author{{\bf Hai-Long Her}\footnote{The work is supported by the Supporting Program for Postdoctoral
Research of Jiangsu Province, P.R.China.}}
\date{}
\maketitle

\begin{quote}
\small {{\bf Abstract}.\  Let $(M,\omega)$ be a compact symplectic
manifold, and $\phi$ be a symplectic diffeomorphism on $M$, we
define a Floer-type homology $FH_*(\phi)$ which is a generalization
of Floer homology for symplectic fixed points defined by Dostoglou
and Salamon for monotone symplectic manifolds. These homology groups
are modules over a suitable Novikov ring and depend only on $\phi$
up to a Hamiltonian isotopy.}
\end{quote}


\section {Introduction.}

Floer homology for symplectic manifolds is a great by-product for
proving the famous Arnold conjectures\cite{A} concerning about the
number of non-degenerate fixed points of  Hamiltonian
diffeomorphisms of any compact symplectic manifold and the number of
transversal intersection points of any Lagrangian submanifold with
its Hamiltonian deformations. For the first case, since the fixed
points of a Hamiltonian diffeomorphism $\psi_H$ correspond to the
time-1 periodic solutions of the Hamiltonian vector field $X_t$
generated by a time-dependent Hamiltonian function $H_t$, the
problem is equivalent to estimating the number of non-degenerate
periodic solutions. To do this, A. Floer\cite{F1}\cite{F2}\cite{F3}
initialed the method of using Gromov's pseudoholomorphic curves as
connecting orbits between the periodic solutions and, by counting
connecting orbits, establishing infinite dimensional Morse-Witten
complex and homology. For monotone (by Floer\cite{F3}) and
semi-positive (by Hofer-Salamon\cite{HS}, Ono\cite{On1}) symplectic
manifolds, defining Floer homology only involves the moduli space of
connecting orbits, and the Gromov's weak compactness for
$J$-holomorphic curves can be proved (with some kind of dimension
counting argument). As for general symplectic manifolds, however,
the phenomenon of bubbling-off spheres can not be avoided by
transversality arguments, and for this reason, the Gromov's
compactness for the moduli space of $J$-holomorphic curves
connecting periodic solutions can not always hold. By considering
the enlarged moduli space of stable maps  and using the elaborate
virtual techniques (or establishing the so-called Kuranishi
structure), Fukaya-Ono\cite{FO}, Liu-Tian\cite{LT} independently
gave the definition of general Floer homology. As a result, the
general Arnold conjecture for nondegerate fixed points of
Hamiltonian diffeomorphism was proved.

By now, there exist some variants of Floer-type homology. For
instance, for (simply connected) monotone symplectic manifold
$(M,\omega)$, using the same method by Floer\cite{F3} with some
modified details, Dostoglou-Salamon\cite{DS} defined a new (or
generalized) version of Floer homology for symplectic fixed points
of some special symplectomorphism as $\phi_H=\psi_1^{-1}\circ\phi$,
where $\psi_1$ is the time-1 map of symplectomorphisms $\psi_t:
M\rightarrow M$ generated by a Hamiltonian $H_t$ satisfying
$H_t=H_{t+1}\circ\phi$, $i.e.$
$$\frac{d}{dt}\psi_t=X_t,\ \ \psi_0=id,\ \ \iota(X_t)\omega=dH_t$$
with $\psi_{t+1}\circ\phi_H=\phi\circ\psi_t$, and $\phi$ is a given
symplectomorphism on $M$. It is easy to see that when $\phi=id$,
$\phi_H$ is a Hamiltonian diffeomorphism, and the generalized Floer
homology is just the ordinary one. For a symplectomorphism $\phi\in
Symp_0(M,\omega)$, $i.e.$ in the identity component of symplectic
diffeomorphism group $Symp(M,\omega)$, Ono\cite{On2} also defined
another version of so-called Floer-Novikov homology for general
symplectic manifold.

In the present paper, we generalize the Dostoglou-Salamon's
construction to general symplectic manifold and define such Floer
homology. The fundamental idea is to study the moduli space of
stable connecting orbits from open surface from which we construct
virtual cycle, and to use it to define a complex which is a module
over a Novikov ring generated by the solutions of the Hamiltonian
equation above, and a suitable boundary operator on it whose
homology is the favorite one.

Firstly, we review the basic idea of dealing with problems of moduli
space.
\begin{Definition}
We say a triple $(\Cal{E},\Cal{B},s)$ is a Fredholm system if the
follwing hold

1) $\Cal{E}$, $\Cal{B}$ are two smooth Banach orbifolds,
$\pi:\Cal{E}\rightarrow\Cal{B}$ is a Banach orbibundle, whose each
fibre $\Cal{E}_x,\ x\in\Cal{B}$ is a Babach space.

2) $s:\Cal{B}\rightarrow\Cal{E}$ is a smooth section, and for
$\forall x\in\Cal{B}$ the linearization of $s$ at $x$
$$D_xs:T_x\Cal{B}\rightarrow\Cal{E}_x$$ is a Fredholm operator and $s^{-1}(0)$ is compact. $s$ is
called a Fredholm section.
\end{Definition}
And the zero set of section $\Cal{M}=s^{-1}(0)$ is called the moduli
space of the system.

In particular, if the linearization map $D_xs$ is surjective for all
$x\in\Cal{M}$, then $\Cal{M}$ is a smooth orbifold and the dimension
of $\Cal{M}$ is equal to the Fredholm index of $D_xs$. If
$ind(D_xs)=d$, then $\Cal M$ can be considered as a cycle in
$H_d(\Cal{B})$ representing the Euler class of the bundle
$\Cal{E}\rightarrow\Cal{B}$. Thus some invariants (for example G-W
invariants), can be defined for any cohomology form $\alpha\in
H^d(\Cal{B},\R)$ as $\Phi(\alpha)=\int_{\Cal{M}}\alpha$.

Actually, defining Floer homology for fixed point of Hamiltonian
diffeomorphism $\Phi_H$ is also based on the study of thus moduli
space. For instance, let $\Cal{B}={\Cal W}^{1,p}(u^*TM)$ be some
suitable completion of the space of smooth maps $u:\Sigma\rightarrow
(M,\omega)$ with some suitable boundary condition, where $\Sigma$ is
a Riemann surface (for example $\Sigma=\R\times S^1$ and
$\lim_{s\rightarrow\pm\infty}u(s,t)=x_{\pm}(t)$ are periodic
solutions of the Hamiltonian equation), and let $\Cal{E}$ be the
Banach bundle over $\Cal{B}$, whose fiber at $u\in\Cal{B}$ is the
Banach space $L^p(\Lambda^{0,1}(u^*TM))$, then the
$\bar{\partial}_J$ map or its some perturbation
$\bar{\partial}_{J,H}$ can be thought as a section of bundle
$\Cal{E}\rightarrow\Cal{B}$. In some ideal case, for instance $M$ is
a monotone symplectic manifold, Floer proved that
$(\Cal{E},\Cal{B},\bar{\partial}_{J,H})$ is a Fredholm system, and
for generic chosen pair $(J,H)$, the linearization operator $D_u$ of
$\bar{\partial}_{J,H}$ is surjective at any $u\in{\Cal
M}=\bar{\partial}_{J,H}^{-1}(0)$. Using this moduli space $\Cal M$,
especially, by studying its 1 and 2 dimensional components, Floer
defined his well known homology.

However, in general, the moduli space $\Cal M$ is not compact.
Floer's method for defining homology for general symplectic manifold
is invalid, since we can not naturally get a Fredholm system
$(\Cal{E},\Cal{B},\bar{\partial}_{J,H})$ as above. In fact, one can
overcome this difficult by considering a larger space of maps, {\it
i.e.} the stable compactification $\overline{\Cal M}$, say moduli
space of stable maps, which is the zero set of an enlarged bundle,
introduced by Kontsevich. At the same time, although the new
enlarged moduli space is compact, the appearing of multiple covered
$J$-spheres with negative Chern number would make the linearization
map $D_xs$ be not always surjective for any $x\in \overline{\Cal
M}$. This also makes the dimension of the ``boundary" of
$\overline{\Cal M}$ too larger than estimated by index theorem. The
idea of virtual techniques to deal with this problem is to construct
generic perturbative section $\bar{\partial}_{J,H}+\nu$ of some
orbibundle (or, say multi-bundle) over the enlarged space of stable
${\Cal W}^{1,p}$-maps, with the virtual moduli space ${\Cal M}^\nu$
as the zero set of this perturbed section, and to construct the
so-called virtual moduli cycle $C({\Cal M}^\nu)$, from which one can
derive the well-defined Floer-type homology.

Now in our setting, for any two $\tilde{x}_-,\ \tilde{x}_+$ the
stable moduli space is $P{\Cal M}(\tilde{x}_-, \tilde{x}_+)$, where
we denote $\tilde{x}_- (\tilde{x}_+)$ for the lift of fixed point
$x_-(x_+)$ of $\phi_H$ in some universal covering space. Roughly
speaking, $P{\Cal M}(\tilde{x}_-, \tilde{x}_+)$ consists of all
stable $(J,H)$-orbits connecting $\tilde{x}_-$ and $\tilde{x}_+$,
such a stable connecting orbit $V: \Sigma\rightarrow M$ contains
some main components $v_m,\ m=1,\cdots,K$, which each is a
$\bar{\partial}_{J,H}$-orbit, and some bubble components $f_b,\
b=1,\cdots,L$, which each is a $J$-holomorphic sphere. And $P{\Cal
M}(\tilde{x}_-, \tilde{x}_+)$ is the natural compactification of the
ordinary moduli space ${\Cal M}(\tilde{x}_-, \tilde{x}_+)$ which
consists of only $\bar{\partial}_{J,H}$-orbits with open domain
$\R^2$. We refer the reader to section 2 and 3 for related
definitions. For monotone symplectic manifold $M$, Dostoglou-Salamon
used only the moduli space ${\Cal M}(\tilde{x}_-, \tilde{x}_+)$ to
define the homology for symplectomorphism $\phi$, while for general
symplectic manifold we must consider its stable compactification.

Then the ambient space is denoted by ${\Cal B}(\tilde{x}_-,
\tilde{x}_+)$, and in the partially smooth category (c.f. \cite{M99}
or section 4) we show that there exists a neighborhood $\Cal W$ of
$P{\Cal M}(\tilde{x}_-, \tilde{x}_+)$ in ${\Cal B}(\tilde{x}_-,
\tilde{x}_+)$ with the so-called multi-fold structure or atlas
$\widetilde{\Cal V}$, and we define a multi-bundle $\widetilde{\Cal
E}$ over it, then we can show that in a small covering neighborhood
${\Cal W}_\epsilon$, locally the generic perturbed map
$\bar{\partial}_{J,H}+\nu$ gives the ``fine" section, and all of
them fit together to give a multi-section $\tilde{s}$ of this
bundle, In other words, for generic pair $(J,H)$ we can obtain a
transverse Fredholm system $(\widetilde{\Cal E},\widetilde{\Cal
V},\Cal{W},\tilde{s})$ with index $d$, and the zero set of this
multi-section is the virtual moduli space, denoted by $P{\Cal
M}^\nu(\tilde{x}_-, \tilde{x}_+)$ which is compact and can be
considered as a relative cycle with correct dimension $d$ estimated
by index theorem.

We then define a graded $\Q$-space
$C_*=C_*(J,H,\phi)=\oplus_nC_n(J,H,\phi)$  as usual. Simply
speaking, we can define a functional $F$ on some covering of a path
space, and $C_n(J,H,\phi)$ is generated by the critical points of
$F$ with Conley-Zehnder index $n$. And $C_*$ is in general an
infinite dimensional space over $\Q$ but a finite dimensional space
over a Novikov ring $\Lambda_{\omega,\phi}$ (c.f. section 6.2 for
details). Then we can just define the boundary operator by
$\delta_{J,H,\nu}:\ C_n\rightarrow C_{n-1}$, such that for any
$\tilde{x}\in Crit_n(F)$,
$$\delta_{J,H,\nu}(\tilde{x})=\sum_{\tilde{y}\in Crit_{n-1}(F)}\#(P{\Cal
M}^\nu(\tilde{y},\tilde{x}))\tilde{y},$$ where $Crit_n(F)$ denotes
the set of critical points of $F$ with Conley-Zehnder index $n$. Now
we can state our main result
\begin{Theorem}\label{mainthm}
Let $(M,\omega)$ be a compact symplectic manifold with compatible
almost complex structure $J$, and $\phi$ be a given
symplectomorphism on $M$. Then for a generic pair of $J$ and time
dependent Hamiltonion function $H:\R\times M\rightarrow\R$
satisfying $H_t=H_{t+1}\circ\phi$, we can construct a compact
relative rational cycle $P{\Cal M}^\nu$ and a boundary operator
$\delta_\nu=\delta_{J,H,\nu}:C_*\rightarrow C_*$, such that
$\delta_\nu^2=0$. So the Floer homology $FH_*(J,H,\phi,\nu)$ is the
homology of this chain complex $(C_*,\delta_\nu)$. Moreover, Floer
homology $FH_*(J,H,\phi,\nu)$ is independent of the choice of the
generic pair $(J,H)$, and it depends on the symplectomorphism $\phi$
only up to Hamiltonian isotopy, i.e. there exists a natural
isomorphism $$FH_*(J_0,H_0,\phi_0,\nu_0)\rightarrow
FH_*(J_1,H_1,\phi_1,\nu_1),$$ provided $\phi_0$ and $\phi_1$ are
related by a Hamiltonian isotopy.
\end{Theorem}
{\it Remark}. With the Floer homology $FH_*(\phi)$ defined above, we
can consider the {\it pair-of-pants} construction suggested by
Donaldson, {\it i.e.} we can consider a homomorphism
$$FH_*(\psi)\otimes FH_*(\phi)\rightarrow FH_*(\psi\circ\phi).$$ If
$\psi=id$, this induces the quantum cap product. We will not study
this topic in the present article, the author plans to treat the
quantum cap product and its related applications in another paper.

In the present paper, we also use the similar virtual techniques
with some modifications based on Liu-Tian's work \cite{LT} to
construct Floer homology of symplectomorphism for general symplectic
manifold. In fact, the methods by Fukaya-Ono\cite{FO} and Li-Tian
\cite{LiT} can be also used for this purpose, however, we want to
use as less analytic tools as possible at the cost of more
topological arguments. Also the technique used by Ruan \cite{R} may
be applied for this by generalizing arguments in \cite{Si} to
construct virtual neighborhoods.

The contents of the paper are as follows. In section 2 we define the
ordinary moduli space ${\Cal M}(\tilde{x}_-,\tilde{x}_+)$ of
$\bar{\partial}_{J,H}$-connecting orbits with the variant in that
the solutions of the Hamiltonian equation given above, whose lifts
are considered as the critical points of action functional, are
unnecessary periodic ones. In section 3 we define the stable
connecting orbits and the enlarged stable moduli space $P{\Cal
M}(\tilde{x}_-,\tilde{x}_+)$. Then in section 4 we show that this
enlarged space is Hausdorff and compact, and we study its small
neighborhood in its ambient space. In section 5, in a more abstract
setting and in the partially smooth category, we show the concepts
of multi-fold, multi-bundle and multi-section which can be found
(maybe with slight difference in different settings) in many
literatures  (cf. \cite{FO}\cite{LT}\cite{M99}\cite{S}\cite{Si}),
and show the general method of constructing virtual cycle. In
section 6 we prove that locally we can get a transverse Fredholm
system by generically perturbation. As an application of the results
in preceding sections, in section 7 we construct the virtual moduli
space and the relative virtual moduli cycle $P{\Cal M}^\nu$ and
define the Floer-type homology.

\bigskip
\noindent {\bf Acknowledgement}. The author is grateful to Professor
Gang Tian for many advices and constant encouragement in his work.

\section {Moduli space of connecting orbits.}

We take the twisted free loop space as
$$\Omega_{\phi}=\{\gamma:
\R\rightarrow M: \gamma(t+1)=\phi(\gamma(t))\}.$$ And consider the
closed 1-form on $\Omega_{\phi}$
$$\langle\alpha(\gamma),\xi\rangle=\int^1_0\omega(\dot{\gamma}-X_t(\gamma),\xi)dt.$$

We say that an almost complex structure $J:TM\rightarrow TM$ is
compatible with symplectic form $\omega$ if they induce a Riemannian
metric $g(\xi,\eta)=\omega(\xi,J\eta)$. And we denote the set of all
compatible almost complex structure by ${\Cal J}(M,\omega)$. Choose
a smooth map $\R\rightarrow {\Cal J}(M,\omega): \ t\mapsto J_t$ such
that $J_t=\phi^*J_{t+1}$. Such a structure determines a metric on
$\Omega_{\phi}$.

Now, we consider the minimal covering $\pi:\
\tilde{\Omega}_{\phi}\rightarrow\Omega_{\phi}$ such that the form
$\pi^*\alpha$ is exact, $i.e.$ there is a functional $F$ on
$\tilde{\Omega}_{\phi}$, such that $\ \pi^*\alpha=dF$, and its
structure group $\Gamma$ is free abelian. In general, we may
additionally assume that there is an injective homomorphism
$i:\Gamma \rightarrow\pi_2(M)$. For instance, if the symplectic
manifold $M$ is simply connected, then
$\Gamma=\pi_1(\Omega_\phi)\cong \pi_2(M)$.

In this article, we will at first consider a connected component of
$\Omega_{\phi}$, and describe a certain covering space of it.
Firstly, we choose and fix a path $\gamma_0\in\Omega_{\phi}$, and
consider the component of $\Omega_{\phi}$ containing $\gamma_0$,
denote it by $\Omega_{\phi}^{\gamma_0}$ or $\Omega_{\phi}^{0}$. We
denote by $C([0,1]\times\R,M)$  for all the continuous maps
$w:[0,1]\times\R\rightarrow M$, we consider the set of pairs as
$$\{(\gamma,w)|\ w(0,\cdot)=\gamma_0,\ w(1,\cdot)=\gamma,\ \gamma\in\ \Omega_\phi^0,\ w\in C([0,1]\times\R,M)\}.$$
We define a homomorphism $\Cal {I}_\omega: \pi_{2} (M)\rightarrow
\R$, ${\Cal I}_\omega(A)=\int_A\omega$, $\forall A\in\pi_2(M)$. And
we define a weaker equivalence relation by
$$(\gamma,w)\sim(\gamma',w')\Longleftrightarrow\gamma=\gamma'\
{\rm and}\ \int_{[0,1]\times[0,1]} w^*\omega=\int_{[0,1]\times[0,1]} (w')^*\omega.$$

The universal covering space of $\Omega_\phi^0$ can be defined as
the set of equivalence classes of the pairs defined above, denote it
by $\tilde{\Omega}^0_{\phi}=\{[\gamma,w]\}$. Then the $\Gamma$ is
the group of deck transformation of the covering
$$\tilde{\Omega}^0_{\phi}\rightarrow\Omega_\phi^0.$$ One can check
that this is a regular covering, say,
$\tilde{\Omega}^0_{\phi}/\Omega_\phi^0\simeq\Gamma$ for any fixed
$\gamma_0\in\Omega_\phi^0$. Furthermore, we see that $\pi_2(M)$ acts
on $\tilde{\Omega}_{\phi}$ via
$\tilde{\gamma}=[\gamma,w]\rightarrow[\gamma,A\#w]=\tilde{\gamma}\#A$,
where $A\#w$ denotes the equivalence class of the connected sum
$a\#w$ for any representative $a$ of $A$. The covering map is
$\pi:\tilde{\Omega}_{\phi}^0\rightarrow \Omega_{\phi}^0$,
$\pi(\tilde{\gamma})=\gamma$. Now, we can write the functional
$F_{\gamma_0}:\tilde{\Omega}^0_{\phi}\rightarrow \R$ by
$$F_{\gamma_0}([\gamma,w])=\int_{[0,1]\times[0,1]} w^*\omega-\int_0^1H(\gamma(t),t) dt.$$
Then the universal cover $\tilde{\Omega}_{\phi}$ of $\Omega_{\phi}$,
and the functional $F$ on $\tilde{\Omega}_{\phi}$ can be defined
componentwise. The critical points of $F$ are such $[\gamma,w]$ with
$\gamma$ being the solution of Hamiltonian equation listed in the
introduction.\footnote{For a special case $\phi=id$, our covering
space and functional $F$ are different from the usual covering space
$\widetilde{L}$ and action functional $a_H$ ({e.g.
\cite{HS}\cite{LT}}), however, in the path space $\Omega_{\phi}$
they have the same paths as the projection of the set of critical
points. The author thinks it is interesting to consider the relation
between the two kinds of Floer homology constructed under the two
settings.}

In the sequel, for simplicity of notations, we will just write
$\Omega_{\phi}$ and $\tilde{\Omega}_{\phi}$ for $\Omega_{\phi}^0$
and $\tilde{\Omega}_{\phi}^0$, respectively, if without the danger
of confusion. We also use $\tilde{x}$ to denote a critical point of
the functional $F$, which can be considered as a lift of the
symplectic path $\psi_t(x)$ or a lift of a fixed point $x$ of the
symplectomorphism $\phi_H$.

The gradient $\nabla F$ of the functional $F$, with respect to the
lift of the metric on $\Omega_{\phi}$, is a $\Gamma$-invariant
vector field on $\tilde{\Omega}_{\phi}$, and $\pi_*\nabla
F=grad_{\alpha}$. Then we consider the moduli space of thus gradient
flows connecting a pair of critical points
$(\tilde{x}_-,\tilde{x}_+)$ of $F$
$$\widehat{\Cal M}(\tilde{x}_-,\tilde{x}_+)=\left\{
\begin{array}{ccc}
\tilde{u}:\ \R\rightarrow\tilde{\Omega}| \ &
\frac{d\tilde{u}(s)}{ds}=-\nabla F(\tilde{u}(s)),\ \tilde{u}\ {\rm
is\ not\ constant},\\

 & \lim_{s\rightarrow\pm\infty}\tilde{u}(s)=\tilde{x}_{\pm}
\end{array}
\right\}.
$$

It is also called the moduli space of connecting orbits. Denote the
collection by $\widehat{\Cal
M}=\bigcup_{\tilde{x}_{\pm}}\widehat{\Cal
M}(\tilde{x}_-,\tilde{x}_+),$ the non-parameterized space by $${\Cal
M}(\tilde{x}_-,\tilde{x}_+)=\widehat{\Cal
M}(\tilde{x}_-,\tilde{x}_+)/\R,$$ and the natural quotient map by
$q:\ \widehat{\Cal M}\rightarrow{\Cal M} $. For generic pair
$(H,J_t)$, all ${\Cal M}(\tilde{x}_-,\tilde{x}_+)$ are finite
dimensional smooth manifolds.

The map $u:\ \R^2\rightarrow M$, defined by
$u(s,t)=\pi(\tilde{u}(s))(t)$, satisfies the following perturbed
Cauchy-Riemann equation
\begin{equation}\label{001}
\bar{\partial}_{J,H}(u)=\partial_s u+J_t(u)(\partial_tu-X_t(u))=0
\end{equation}
with boundary condition
\begin{equation}\label{002}
u(s,t+1)=\phi(u(s,t)),
\end{equation}
and has limits
\begin{equation}\label{003}
\lim_{s\rightarrow\pm\infty}u(s,t)=\psi_t(x_{\pm}),\ \
x_{\pm}=\phi_H(x_{\pm}),
\end{equation}
where $x_{\pm}$ are nondegenerate fixed points of $\phi_H$.

The energy of $u$ is
$$E(u)=\frac{1}{2}\int_{-\infty}^{\infty}\int_0^1(|\partial_su|^2+|\partial_tu-X_t(u)|^2)dtds<\infty$$
It is easy to see that
$\tilde{x}_-\#u=[x_+,w\#u]=[x_+,w']=\tilde{x}_+$, where $w'=w\#u$ is
the obvious concatenation of $w$ and $u$ along the path $x_-$. And
we have the equality
\begin{equation}\label{004}
E(u)=F(\tilde{x}_+)-F(\tilde{x}_-).
\end{equation}

For any smooth function $u:\R^2\rightarrow M$ satisfying
$u(s,t+1)=\phi(u(s,t))$, we denote $W_{\phi}^{k,p}(u^*TM)$ for the
Sobolev completion of the space of smooth  vector fields
$\xi(s,t)\in T_{u(s,t)}M$ along $u$, which satisfy
$\xi(s,t+1)=d\phi(u(s,t))\xi(s,t)$ and have compact support on
$\R\times [0,1]$, with respect to the $W^{k,p}$-norm over
$\R\times[0,1]$. And denote
$L_{\phi}^{p}(u^*TM)=W_{\phi}^{0,p}(u^*TM)$. Then for solution $u$
of (\ref{001}) and (\ref{002}), we get the following linear operator
by linearizing (\ref{001})
$$D_u: W_{\phi}^{1,p}(u^*TM)\rightarrow L^{p}_{\phi}(u^*TM)$$ defined
by
$$D_u\xi=\nabla_s\xi+J_t(u)(\nabla_t\xi-\nabla_{\xi}X_t(u))+\nabla_{\xi}J_t(u)(\partial_tu-X_t(u)),$$
where $\nabla$ denotes the Levi-Civita connection associated to the
metric $g_t(\cdot,\cdot)=\omega(\cdot,J_t\cdot)$. If the fixed
points $x_{\pm}=\phi_H(x_{\pm})$ are nondegenerate then $D_u$ is a
Fredholm operator and its index is given by the Maslov index of $u$.

More precisely, we follow \cite{DS1} to show the definition of the
Maslov index. Let $Sp(2n)$ be the group of symplectic matrices. And
denote the singular subset by $Sing(2n)=\{A\in Sp(2n)|\ {\rm
det}(Id-A)=0\}$, which is called {\it Maslov cycle}, its complement
is an open and dense subset of $Sp(2n)$, denote it by
$$Sp^*(2n)=Sp(2n)-Sing(2n).$$ For any path $\Psi: [0,1]\rightarrow
Sp(2n)$, with $\Psi(0)=Id$ and $\Psi(1)=Sp^*(2n)$, there exists the
so-called Conley-Zehnder index $\mu_{CZ}(\Psi)$ (c.f.
\cite{RS}\cite{SZ}). Given two nondegenerate fixed points $x_{\pm}$
of $\phi_{H}$, denote by ${\Cal P}(x_-,x_+)$ the space of all smooth
functions $u:\R^2\rightarrow M$ satisfying (\ref{002}) and
(\ref{003}). For any $u\in {\Cal P}(x_-,x_+)$, we choose a
trivialization $\Phi(s,t):\R^{2n}\rightarrow T_{u(s,t)}M$ such that
$$\Phi(s,t)^*\omega=\omega_0,\ \ \ \
\Phi(s,t+1)=d\phi(u(s,t))\Phi(s,t)$$and
$\lim_{s\rightarrow\pm\infty}\Phi(s,t)=\Phi^{\pm}(t):\
\R^{2n}\rightarrow T_{\psi_t(x_{\pm})}M$. Now, we construct the
symplectic paths
$$\Psi^{\pm}(t)=\Phi^{\pm}(t)^{-1}d\psi_t(x_{\pm})\Phi^{\pm}(0).$$
It is easy to see that $\Psi^{\pm}(1)\in Sp^*(2n)$ since $x_{\pm}$
are nondegenerate. Then we can define the Conley-Zehnder index of
$x_\pm$  by
$$\mu_{CZ}(x_\pm)=\mu_{CZ}(\Psi^\pm),$$ and the Maslov index of the
pair $(u,H)$ is defined by
$$\mu(u,H)=\mu_{CZ}(x_+)-\mu_{CZ}(x_-)=\mu_{CZ}(\Psi^+)-\mu_{CZ}(\Psi^-),$$ which is independent
of the choice of the trivialization and satisfies the
following properties\\

\noindent ({\it Homotopy}): For given $H$ and two nondegenerate
fixed points $x_{\pm}$ of $\phi_H$, $\mu(u,H)$ is constant on the
homotopy components of ${\Cal P}(x_-,x_+)$.

\noindent ({\it Zero}): If $x_-=x_+$, then $\mu(u,H)=0$.

\noindent ({\it Catenation}):
$\mu(u_{01}\#u_{12},H)=\mu(u_{01},H)+\mu(u_{12},H)$.

\noindent ({\it Chern class}): For $v: S^2\rightarrow M$,
$\mu(u\#v,H)=\mu(u,H)-2c_1(v)$.

\noindent ({\it Morse index}): If $\phi=Id$ and $H_t=H$ is a Morse
function with sufficiently small second derivatives, then fixed
points are the critical points of $H$ and
$\mu(u,H)=Ind_H(x_{+})-Ind_H(x_{-})$.

\noindent ({\it Fixed point index}): For $u\in {\Cal P}(x_-,x_+)$,
$$(-1)^{\mu(u,H)}={\rm sign\ of\ det}(Id-d\phi(x_+)){\rm
det}(Id-d\phi(x_-)).$$

Let $\widetilde{{\rm Fix}}(\phi_H)\subset\widetilde{\Omega}_\phi$
denote the set of elements which cover curves of the form
$x(t)=\psi_t(x)$, $x\in{\rm Fix}(\phi_H)$. Then we have a fibration
of discrete sets
$$\Gamma\hookrightarrow\widetilde{{\rm Fix}}(\phi_H)\rightarrow{\rm
Fix}(\phi_H).$$ Every function $u\in\Cal{P}(x_-,x_+)$ and every lift
$\tilde{x}_-\in\widetilde{{\rm Fix}}(\phi_H)$ of $x_-$ determines a
unique lift $\tilde{x}_+=\tilde{x}_-\#u$ of $x_+$. By the
``homotopy" and ``catenetion" properties of the Maslov index there
exists a unique map $\mu_{rel}:\widetilde{{\rm Fix}}(\phi_H)\times
\widetilde{{\rm Fix}}(\phi_H)\rightarrow \Z$ such that
$$\mu(u,H)=\mu_{rel}(\tilde{x}_-,\tilde{x}_+)$$ whenever
$\tilde{x}_+=\tilde{x}_-\#u$. Then by the ``Chern class" property
one has
$$\mu_{rel}(\tilde{x}_-,v\#\tilde{x}_+)=\mu_{rel}(\tilde{x}_-,\tilde{x}_+)-2c_1(v)$$
for $v\in\pi_2(M)$, and the ``catenetion" property reads
$$\mu_{rel}(\tilde{x}_0,\tilde{x}_1)+\mu_{rel}(\tilde{x}_1,\tilde{x}_2)=\mu_{rel}(\tilde{x}_0,\tilde{x}_2).$$
Now, if we fix a reference critical point of the functional $F$, say
$\tilde{x}_0\in\widetilde{{\rm Fix}}(\phi_H)$, such that its
projection $x_0\in{\rm Fix}(\phi_H)$ is nondegenerate, then we can
define the Conley-Zehnder index of any $\tilde{x}\in\ Crit(F)$ by
$$\mu_{CZ}(\tilde{x})=\mu_{rel}(\tilde{x}_0,\tilde{x}).$$
It is easy to see that the dimension of the nonparameterized moduli
space ${\Cal M}(\tilde{x}_-,\tilde{x}_+)$ is
$\mu_{rel}(\tilde{x}_-,\tilde{x}_+)-1$.

\section {Moduli space of stable connecting orbits.}

Now we follow the method by Fukaya-Ono\cite{FO} with some
modifications to define the stable connecting orbits and moduli
space. We denote the set of critical points of the function $F$ by
{\it Crit(F)}. For convenience of the reader, we list some
necessary definitions.
\begin{Definition}
A semi-stable curve with k marked points is a pair
$(\Sigma,\bf{z})$, where the set
$\Sigma=\cup\pi_{\Sigma_{\nu}}(\Sigma_{\nu})$ is connected, each
$\Sigma_{\nu}$ is a Riemann surface, and the number of $\Sigma_\nu$
is finite, $\pi_{\Sigma_{\nu}}:\Sigma_{\nu}\rightarrow\Sigma$ is a
continuous and locally homeomorphic map, and ${\bf
z}=(z_1,\cdots,z_k)$ are k distinguished points on $\Sigma$.
Moreover, the following hold

$1^{\circ}$ For each $p\in\Sigma$, ${\rm
Sum}_{\nu}\#\pi_{\Sigma_{\nu}}^{-1}(p)\le 2$; For each $z_i$, ${\rm
Sum}_{\nu}\#\pi_{\Sigma_{\nu}}^{-1}(z_i)=1$.

$2^{\circ}$ The set $\{p|{\rm
Sum}_{\nu}\#\pi_{\Sigma_{\nu}}^{-1}(p)= 2\}$ is of finite order.
\end{Definition}
If ${\rm Sum}_{\nu}\#\pi_{\Sigma_{\nu}}^{-1}(\pi_{\Sigma_{\nu}}(p))=
2$, $p$ is called a {\it singular or double} point. If
$\pi_{\Sigma_{\nu}}(p)=z_i$ for some $i$, we say the
$p\in\Sigma_{\nu}$ is {\it marked}. And say $\Sigma_{\nu}$ is a {\it
component} of $\Sigma$. If all components are spheres, we say
$(\Sigma,{\bf z})$ is a genus 0 semi-stable curve.

A homeomorphism $h:(\Sigma,{\bf z})\rightarrow(\Sigma,{\bf z})$ is
called an {\it automorphism} if it can be lifted to bi-holomorphic
isomorphisms $h_{\mu\nu}:\Sigma_\mu\rightarrow\Sigma_\nu$ for each
component and $h(z_i)=z_i$ for each $i$. We denote the automorphism
group of $(\Sigma,{\bf z})$ by $Aut(\Sigma,{\bf z})$ or $G_\Sigma$.
\begin{Definition}
Let $(M,\omega)$ be a symplectic manifold with a compatible almost
complex structure $J:TM\rightarrow TM$. A continuous map
$u:\Sigma\rightarrow M$ is called {\rm pseudo-holomorphic}  if the
composition $u\circ\pi_{\Sigma_\nu}:\Sigma_\nu\rightarrow M$ is
pseudo-holomorphic (or J-holomorphic ) for each $\nu$.
\end{Definition}

\begin{Definition}\label{3.3}
A pair $((\Sigma,{\bf z}),u)$ is called a $J$-{\rm stable map} if
for each $\nu$ one of the following conditions holds

$1^{\circ}$ $u\circ\pi_{\Sigma_\nu}:\Sigma_\nu\rightarrow M$ is a
nonconstant $J$-holomorphic map.

$2^{\circ}$ Let $m_\nu$ be the number of special points on
$\Sigma_\nu$ which are singular or marked, then $m_\nu\geq 3$.
\end{Definition}
If a semi-stable curve $(\Sigma,\bf{z})$ satisfies the condition
$2^\circ$ in the Definition \ref{3.3} for each component, then we
say it is a {\it stable curve}. For a pair $((\Sigma,{\bf z}),u)$,
we define its {\it automorphism group} by
$$Aut((\Sigma,{\bf z}),u)=\{h:\Sigma\rightarrow\Sigma|\ h\ {\rm is\
an\ automorphism,\ and}\ u\circ h=u \}.$$ It can be proved that
$((\Sigma,{\bf z}),u)$ is stable if and only if $Aut((\Sigma,{\bf
z}),u)$ is a finite group\cite{FO}.

\begin{Definition}
Let $\beta\in H_2(M,\Z)$. We denote by $(M,J,\beta)_{0,k}$ for the
set of genus 0 stable maps $((\Sigma,{\bf z}),u)$ with $k$ marked
points such that $u_*([\Sigma])=\beta$. We say two stable maps are
equivalent, {\it i.e.} $((\Sigma,{\bf z}),u)\sim((\Sigma,{\bf
z}'),u')$, if and only if there exists an isomorphism
$h:(\Sigma,{\bf z})\rightarrow(\Sigma,{\bf z}')$ satisfying $u'\circ
h=u$ and $h(z_i)=z_i'$ for each $i$. We write $C{\Cal
M}_{0,k}(M,J,\beta)=(M,J,\beta)_{0,k}/\sim$, and call it the moduli
space of this kind of stable maps.
\end{Definition}

We define the energy of a genus $0$ stable map $((\Sigma,{\bf
z}),u)\in (M,J,\beta)_{0,k}$ by
$$E(u)=E((\Sigma,{\bf z}),u)=\int_\Sigma u^*\omega$$

Recall that ${\Cal M}(\tilde{x}_-,\tilde{x}_+)$ is the moduli space
of gradient flow lines connecting two critical points
$\tilde{x}_-,\tilde{x}_+$ of $F$. Note the periodicity condition
(\ref{002}), we give the following
\begin{Definition}\label{stable}
A stable connecting orbit between $\tilde{x}_-=\tilde{x}_1$ and
$\tilde{x}_+=\tilde{x}_{K+1}$ is a triple
$((v_1,\cdots,v_K),(f_1,\cdots,f_l),o)$ such that\\

\noindent (1) $v_j=u_j|_{\R\times[0,1]}$, where
$u_j=\pi(\tilde{u}_j)$, $\tilde{u}_j\in {\Cal
M}(\tilde{x}_j,\tilde{x}_{j+1})$, $\tilde{x}_j\in Crit(F)$,
$j=1,\cdots,K+1$.\\

\noindent (2) $f_i=(\Sigma_{f_i}, u_{f_i})\in C{\Cal
M}_{0,1}(M,J,\beta_i)$, where $\Sigma_{f_i}$ is a genus zero
semi-stable curve with one marked point and
$u_{f_i}:\Sigma_{f_i}\rightarrow M$, and
$[u_{f_i}(\Sigma_{f_i})]=\beta_i$. Let $z_i\in\Sigma_{f_i}$ be the marked point.\\

\noindent (3) $o$ is an injection from $\{1,\cdots, l\}$ to the $K$
copies of $\R\times[0,1]$. If $o(i)=(s_i,t_i)$ is on the $j^{\rm
th}$ copy of $\R\times[0,1]$, we require that
$u_{f_i}(z_i)=v_j(s_i,t_i)$. Moreover, if there exists some
$i\in\{1,\cdots,l\}$ satisfying
$o(i)=(s_i,0)$ or $o(i)=(s_i,1)$, then there exists a $j\in\{1,\cdots,l\}$ such that $o(j)=(s_j=s_i,1)$
or $o(j)=(s_j=s_i,0)$, respectively.\\

\noindent (4) If $\tilde{x}_j=\tilde{x}_{j+1}$, then there exists an
i such that $o(i)$ is on the $j^{\rm th}$ copy of $\R\times[0,1]$.
\end{Definition}

We say each $v_j$ is a {\it main component}  and each component of
$f_i$ is a {\it bubble component} of the stable connecting orbit. In
particular, sometimes if necessary, we denote by $f^B_i$ for the
distinct $f_i$ above satisfying $o(i)=(s_i,0)$ or $o(i)=(s_i,1)$ and
by $f^I_i$ for others.

Since each $v_j$ is the restriction of $u_j$, for a sequences
$\{v_j^{(n)}\}$, if there exists a $v_j^*$ such that
$v_j^{(n)}\rightarrow v_j^*$ if and only if there exists a map
$u_j^*$ such that $v_j^*=u_j^*|_{\R\times[0,1]}$ and
$u_j^{(n)}\rightarrow u_j^*$.

Naturally, we define the {\it domain} $\Sigma$ of the stable
connecting orbit as $K$ copies of $\R\times[0,1]$  with two
particular segments of line $L_{\pm\infty}=\{\pm
\infty\}\times[0,1]$, which are called the {\it main components},
denoted by $\Sigma_{m,i},\ i=1,\cdots,K$, attached with $l$ genus
zero semi-stable curves $\Sigma_{f_1},\cdots,\Sigma_{f_l}$ at $l$
points $o(1),\cdots,o(l)$, whose components are called the {\it
bubble components} which each one can be identified with $S^2$, we
assume there totally are $L$ such bubble components, denoted them by
$\Sigma_{b,j},\ j=1,\cdots,L$. And the double points include all
$o(i)$ and the singular points on all $\Sigma_{f_i}$ which are the
intersections of bubble components in $\Sigma_{f_i}$. To make the
domain be stable curve, we can add some marked points on those
semi-stable bubble components such that the condition $2^\circ$ in
Definition \ref{3.3} holds, also in order to consider the
bubbling-off process we maybe add some markings on main components
as well. We call such curve is ${\Cal P}$-stable.

We can simply write the stable connecting orbit as $(\Sigma,V)$,
where $V:\Sigma\rightarrow M$ satisfying $V=v_j$ on the $j^{\rm th}$
copy of $\R\times[0,1]$, and $V=u_{f_i}$ on $\Sigma_{f_i}$.

For ${\bf r}=(r_1,\cdots,r_K)\in\R^K$, it is easy to see that there
exists a natural $\R^K$-action on $(\Sigma,V)$. If $o(i)$ is on the
$j^{\rm th}$ component, also the action is to translate $o(i)$ by
$r_j\in\R$. We say ${\bf r}(\Sigma,V)$ is equivalent to
$(\Sigma,V)$, and denote $[\Sigma,V]$ for the equivalent class of
 $(\Sigma,V)$.

Now we define the moduli space of stable connecting orbits $P{\Cal
M}(\tilde{x}_-,\tilde{x}_+)$,
\begin{Definition}
we say $[\Sigma,V]\in P{\Cal M}(\tilde{x}_-,\tilde{x}_+)$ if
$\tilde{x}_{K+1}=\tilde{x}_+$ and
$$\tilde{x}_1\#(\beta_1+\cdots+\beta_l)=\tilde{x}_-.$$
\end{Definition}

We define the energy of $[\Sigma,V]\in P{\Cal
M}(\tilde{x}_-,\tilde{x}_+)$ is
\begin{equation}\label{005}
E([\Sigma,V])=\sum_{j=1}^K E(u_j)+\sum_{i=1}^{l_1}
E(f^I_i)+\frac{1}{2}\sum_{i=1}^{l_2} E(f^B_i),
\end{equation}
where $l_1=\#\{f^I_i\}$, $l_2=\#\{f^B_i\}$, and $l_1+l_2=l$.

We can see that the energy only depends on
$\tilde{x}_-,\tilde{x}_+,\ i.e.$
$E([\Sigma,V])=F(\tilde{x}_+)-F(\tilde{x}_-)$, and is independent of
the choice of $[(v_1,\cdots,v_K),(f_1,\cdots,f_l),o]\in P{\Cal
M}(\tilde{x}_-,\tilde{x}_+)$.

The following lemma shows that the energy of each genus 0 stable map
or each nonconstant connecting orbit is bounded from below by a
constant.
\begin{Lemma}\label{bound}
Given a compact symplectic manifold $(M,\omega)$ with compatible
almost complex structure $J$. There exists a constant $\delta>0$
such that if $\tilde{x}_-\neq\tilde{x}_+$, then for any nonconstant
genus zero pseudo-holomorphic map $f:\Sigma\rightarrow M$ and
$\forall\ u\in\Cal{M}(\tilde{x}_-,\tilde{x}_+)$,
$\min(E(f),E(u))>\delta$.
\end{Lemma}
The proof is standard, we refer the reader to \cite{G}\cite{FO}.

In order to construct a Fredholm system defined in the Definition
1.1, we define an ambient space of $P{\Cal
M}(\tilde{x}_-,\tilde{x}_+)$. Firstly, we introduce the notion of
stable $W^{k,p}$ or $L^{k,p}$-connecting orbits with
$k-\frac{2}{p}>1$. Given a set $\Sigma$ which is the domain defined
above.
\begin{Definition}\label{W^kp}
We say a triple $((v_1,\cdots,v_K),(f_1,\cdots,f_l),o)$ is a stable
$W^{k,p}$-connecting orbit between $\tilde{x}_-=\tilde{x}_1$ and
$\tilde{x}_+=\tilde{x}_{K+1}$ if

\noindent $1^\circ$ $v_j=u_j|_{\R\times[0,1]}$, where
$u_j=\pi(\tilde{u}_j)$ is a $W^{k,p}$-map such that
$\lim_{s\rightarrow-\infty}\tilde{u}_j(s,\cdot)=\tilde{x}_j$,
$\lim_{s\rightarrow+\infty}\tilde{u}_j(s,\cdot)=\tilde{x}_{j+1}$,
$\tilde{x}_j\in Crit(F)$, $j=1,\cdots,K+1$.

\noindent $2^\circ$ $f_i=(\Sigma_{f_i}, u_{f_i})$, each
$u_{f_i}:\Sigma_{f_i}\rightarrow M$ is a stable $W^{k,p}$-map, and
$[u_{f_i}(\Sigma_{f_i})]=\beta_i$.

\noindent $3^\circ$ The conditions (3) and (4) in the Definition
\ref{stable} are satisfied.

\noindent $4^\circ$ $\tilde{x}_{K+1}=\tilde{x}_+$ and
$\tilde{x}_1\#(\beta_1+\cdots+\beta_l)=\tilde{x}_-$.

\noindent $5^\circ$  Each $W^{k,p}$-map $u_j$ satisfies the
following exponential $W^{k,p}$-decay condition along its ends
$x_j=\pi{\tilde{x}_j}$ and $x_{j+1}=\pi(\tilde{x}_{j+1})$
$$\int\int_{\R\times[0,1]}(|\eta_{j-1}^{(m)}|^p+|\eta_{j}^{(m)}|^p)e^{\varepsilon_0|s|}dsdt<\infty,\ \ m=0,1,\cdots,k,$$
where $\eta_j$ is defined by $u_j(s,t)=exp_{x_j}\eta_j(s,t)$ for
sufficiently large $s$, and $\varepsilon_0$ is a fixed small
positive constant.
\end{Definition}
We denote the moduli space of equivalence classes of all such stable
$W^{k,p}$-orbits connecting $\tilde{x}_-,$ and $\tilde{x}_+$ by
$\Cal{B}(\tilde{x}_-,\tilde{x}_+)$. All other notation introduced
for stable connecting orbits above are also applicable to stable
$W^{k,p}$-orbits. Sometimes, we will simply write $P\Cal{M}$ and
$\Cal{B}$ for these spaces defined above with the ends
$\tilde{x}_\pm$ being not expressed explicitly, and write a stable
orbit, $i.e.$ the triple $((v_1,\cdots,v_k),(f_1,\cdots,f_l),o)$ in
$P\Cal{M}$ and $\Cal{B}$ by $(\Sigma,V)$ or $V$.

\section{Gluing.}
In this section, we will define a small neighborhood $\Cal{W}$ of
the stable moduli space $P{\Cal M}(\tilde{x}_-,\tilde{x}_+)$ in its
ambient space $\Cal{B}(\tilde{x}_-,\tilde{x}_+)$, and we can see on
$\Cal{W}$ there exists a naturally defined topology which is
Hausdorff. The moduli space $P{\Cal M}$ is a compact subset of
$\Cal{W}$. And in some suitably abstract settings, we can consider
$\Cal{W}$ as a space with two topologies, say a {\it partially
smooth} space. On this space one can define a so-called {\it
multi-fold} atlas, and there are related {\it multi-bundles} over it
and a family of compatible multi-valued sections, say {\it
multi-section}, which can be used to define the virtual cycle. The
application of the settings in this section will give rise to the
construction of virtual moduli cycle and the definition of Floer
homology in the section 6.

In order to study the moduli space $P{\Cal M}$, we have to give a
description of the domain. Indeed, we will see that there exist
stratifications $P{\Cal M}=\cup_D P{\Cal M}^D$ and $\Cal{W}=\cup_D
\Cal{W}^D$. The problem is to describe how the strata $\Cal{W}^D$
fit together as the topological type of the domain changes. The
gluing method can give the local ``cornered" coordinate chart (or
say uniformizer) of the space of domains.

\subsection{Gluing the domains}
Firstly, we give a short description of the domains of stable
connecting orbits and the structure of moduli space $\overline{\Cal
PM}_{0,k}$ of such open stable curves, which we call
$\Cal{P}$-stable curves.
Following the notion above, we denote such a $k$-pointed
$\Cal{P}$-stable curve by $(\Sigma,{\bf
z})=(\Sigma,z_1,\cdots,z_k)$. Recall that $\Sigma$ is the union of
$K$ main components $\Sigma_{m,i}, \ i=1,\cdots,K,$ which are copies
of $\R\times[0,1]$ and $L$ bubble components $\Sigma_{b,j},\
j=1,\cdots,L$, which are identified with sphere. Such a stable curve
is said to be $\Cal{P}$-stable if the following things hold. Each
$\Sigma_{m,i}$ has two particular segments of line
$L_{i,\pm\infty}\simeq\{\pm \infty\}\times[0,1]$. All main
components together form a chain such that
$L_{i,+\infty}=L_{i+1,-\infty},\ i=1,\cdots,K-1$.
Sometimes we only write $\Sigma_m$ and $\Sigma_b$ for a main and a
bubble component if without danger of confusion.

We say a main component $\Sigma_m$ of $\Sigma$ is free if it has no
double point, and say a bubble component $\Sigma_b$ of $\Sigma$ is
free if it has at most two double points. In order to get a stable
curve $\Sigma^s$ with minimal marked points, we may at first use the
``forgetting marking" procedure to send $(\Sigma,z_1,\cdots,z_k)$ to
$\Sigma^u$ by ignoring all marked points $z_i,i=1,\cdots,k$, then
add one or two markings on each free bubble component of $\Sigma^u$
to make it stable.

We denote by $G_m$, $G_b$ the automorphism group of $\Sigma_m$ and
$\Sigma_b$, respectively. Note that on the free main component,
$G_m$ is the group of all $\R$-translations acting on $\Sigma_m$,
and on the free bubble component, $G_b$ is the holomorphic spherical
automorphism group preserving double points of $\Sigma_b$. Note that
the automorphism group $G_\Sigma$ of $\Sigma$ consists of all
holomorphic isomorphisms of $\Sigma$ after forgetting its marked
points, and it may interchange different components of $\Sigma$, we
call it the reparameterization group. It contains $G_m$ and $G_b$ as
subgroups.

Two $\Cal{P}$-stable curves $\Sigma_1$ and $\Sigma_2$ are called to
be equivalent if there is a homeomorphism
$\varphi:\Sigma_1\rightarrow\Sigma_2$ preserving the marked points
and the boundary lines of all main components such that the
restriction of $\varphi$ to each component of $\Sigma_1$ is a
holomorphic map. We denote by $\Cal{P}M_{0,k}$ the collection of all
equivalence class of $k$-pointed $\Cal{P}$-stable curves with only
one main component and no bubble component. Its stable
compactification is denoted by $\overline{\Cal{P}M}_{0,k}$, which is
just all $\Cal{P}$-stable curves defined above. Roughly speaking, we
can use the splitting process of main components and let some of
marked points go together or go to the boundary of $\R\times [0,1]$
symmetrically to obtain $\overline{\Cal{P}M}_{0,k}$ from
$\Cal{P}M_{0,k}$.

The topological type of the curve $\Sigma$ is determined by its
number of main components and intersection pattern, say
$I_{\Sigma}$, which pairwisely corresponds to the specific lines
$L_{i,+\infty}=L_{i+1,-\infty}$ in main components and points in the
smooth resolution $\tilde{\Sigma}$ of $\Sigma$ that correspond just
those double points in $\Sigma$. There are of course finite many
main components and intersection patterns. We use the notation
$$\Cal{I}=\{I_{\Sigma} | \Sigma\in \overline{\Cal{P}M}_{0,k} \}$$
to describe the collection of all topological types.

$\overline{\Cal{P}M}_{0,k}$ is stratified according to the
topological type $I_{\Sigma}$ of a curve $\Sigma$ in a stratum. We
can write
$$\overline{\Cal{P}M}_{0,k}=\cup_{I\in\Cal{I}}\Cal{P}M^{I}_{0,k},$$
where $\Cal{P}M^{I}_{0,k}$ is the collection of curves with fixed
intersection pattern $I_{\Sigma}=I$. Each $\Cal{P}M^{I}_{0,k}$ is a
smooth manifold.

The element $\Sigma\in \overline{\Cal{P}M}_{0,k}$ will naturally
appear as domains of the stable $(J,H)$-connecting orbits. We can
always start from the case $\Sigma^s$ which is $\Sigma$ equipped
with minimal number of marked points needed for stability. In the
bubbling process, the topological type of domains maybe change, then
we will also consider those stable curves with extra markings.

Now for a fixed  $I=I_{\Sigma^s}$ and a fixed element $\Sigma\in
\Cal{P}M^{I}_{0,k}$, we follow the idea of Liu-Tian to give a local
description of the nearby curves in $\overline{\Cal{P}M}_{0,k}$.

Since $\Sigma=\Sigma^s$, there is no extra marked point. That is to
say, there are at most two marked points on each component
$\Sigma_l$. We denote the double points by $d_{l,k}$. So the
locations of all double points for the nearby $\Sigma'$ can be
regarded as a local coordinate of $\Cal{P}M^{I}_{0,k}$ near
$\Sigma$. If the double point $d_{l,k}$ lies in the inner part of
the component $\Sigma_l$, we just let $\alpha_{l,k}\in
D_{\delta}(d_{l,k})$ be the complex coordinate of the $\delta$-disc
centering at $d_{l,k}$. If otherwise $d_{l,k}$ lies in the boundary
of a main component $\Sigma_l\simeq \R\times [0,1]$ (certainly
$d_{l,k}=d_{l',k'}$ maybe simultaneously lie in the bubble component
$\Sigma_{l'}$), then let $\alpha_{l,k}$ be the complex coordinate of
the closed half $\delta$-disc centering at $d_{l,k}$, denoted by
$HD_\delta(d_{l,k})$. From the $\varphi$-periodicity condition
(\ref{002}) for the solution of perturbed Cauchy-Riemann equation
(\ref{001}) we know that if bubbling-off occur in the boundary of a
main component $\Sigma_l$, they should appear synchronously on the
two sides of $\Sigma_l\simeq \R\times [0,1]$ satisfying also a
$\varphi$-periodicity. So if $d_{l,k }=(s,\{0\})$ naturally we have
the other double point $d_{l,k'}=(s,\{1\})$ with complex coordinate
$\alpha_{l,k'}$ in the other closed half $\delta'$-disc centering at
$d_{l,k'}$. The collection $\alpha=(\alpha_{l,k})$ is the local
coordinate of $\Cal{P}M^{I}_{0,k}$ near $\Sigma$. The corresponding
curve is denoted by $\Sigma_{\alpha}$.

Then for each double point on the nearby curve
$\Sigma'=\Sigma_{\alpha}$, corresponding to two intersecting
components $\Sigma'_{l_1}$ and $\Sigma'_{l_2}$, we have
$d'_{l_1,k_1}=d'_{l_2,k_2}$. If one double point, say $d'_{l_1,k_1}$
is on the boundary of a main component, we can just temporarily
consider a little larger domain containing $\Sigma'_{l_1}$, $i.e.$
with a larger component $\bar{\Sigma}'_{l_1}\simeq\R\times
[0-\epsilon,1+\epsilon]$\footnote{We have stated that such double
points will appear simultaneously in pair at the two sides of the
boundary, so the following gluing operating will be done
simultaneously for the other double point near the other side of
boundary of $\Sigma'_{l_1}$. }. Then we can set a complex gluing
parameter $t_{l_1,k_1}=t_{l_2,k_2}$ in the $\delta$-disc centering
at the origin of $\C$, and for each pair of lines
$L_{i,+\infty}=L_{i+1,-\infty}$ in two connecting main components
$\Sigma'_i$ and $\Sigma'_{i+1}$, a real gluing parameter
$\tau_i\in[0,\delta]$. Denote the collection of all parameters by
$(\alpha,t,\tau)=\{(\alpha_{l,k},t_{l,k},\tau_i)\}$. Then the
following procedure shows how to get a curve
$\Sigma_{(\alpha,t,\tau)}$ with different topological type from the
curve $\Sigma'=\Sigma_{\alpha}$.

For each double point $d'_{l_1,k_1}=d'_{l_2,k_2}$ of $\Sigma_\alpha$
with coordinates $\alpha_{l_1,k_1}$ and $\alpha_{l_2,k_2}$, take
complex coordinates $z_{l_1,k_1}$ and $z_{l_2,k_2}$ in the two discs
$D_{\delta'}(\alpha_{l_1,k_1})\subset\Sigma'_{l_1}$ (or
$\bar{\Sigma}'_{l_1}$) and
$D_{\delta'}(\alpha_{l_2,k_2})\subset\Sigma'_{l_2}$ respectively.
Suppose $z=e^{-2\pi(r+{\rm i}\theta)}$, then $(r,\theta)$ is the
corresponding cylindrical coordinate. We firstly cut off discs
$$\{(r_{l_1,k_1},\theta_{l_1,k_1})|\ r_{l_1,k_1}>-{\rm
log}|t_{l_1,k_1}|\}\ \ {\rm in}\ \ D_{\delta'}(\alpha_{l_1,k_1})$$
and $$\{(r_{l_2,k_2},\theta_{l_2,k_2})|\ r_{l_2,k_2}>-{\rm
log}|t_{l_2,k_2}|\}\ \ {\rm in}\ \ D_{\delta'}(\alpha_{l_2,k_2}),$$
then along their boundaries and according to the formula
$$\theta_{l_1,k_1}=\theta_{l_2,k_2}+{\rm arg}(t_{l_1,k_1}=t_{l_2,k_2}),$$
we glue back the remaining parts of $D_{\delta'}(\alpha_{l_1,k_1})$
and $D_{\delta'}(\alpha_{l_2,k_2})$. If $d'_{l_1,k_1}$ is on the
boundary of the main component $\Sigma'_{l_1,k_1}$, after above
gluing procedure we only get a mid-step curve
$\bar{\Sigma}_{(\alpha,t)}$ with some new larger main component
$\bar{\Sigma}_{(\alpha,t),l_1}$, one more thing we should do is to
cut off the additional margin of $\bar{\Sigma}_{(\alpha,t),l_1}$ to
obtain a suitable main component $\Sigma_{(\alpha,t),l_1}$.

To use the real parameter $\tau$ to glue two connecting main
components $\Sigma'_i$ and $\Sigma'_{i+1}$ of $\Sigma_\alpha$ along
the line $L_{i,+\infty}=L_{i+1,-\infty}$ is much simpler and direct.
Let $\tau_i^+=\tau_{i+1}^-$, we cut off a strip
$[\frac{1}{\tau_i^+},+\infty)\times[0,1]$ in the main component
$\Sigma_i'$ and a strip
$(-\infty,\frac{1}{\tau_{i+1}^-}]\times[0,1]$ in $\Sigma_{i+1}'$,
then we just simply glue the two remaining parts with identifying
the two end lines.

So the resulting curve is just $\Sigma_{(\alpha,t,\tau)}$, which is
an element of $\overline{\Cal{P}M}_{0,k}$ near $\Sigma$. The
parameter $(\alpha,t,\tau)$ is a ``cornered" coordinate chart of
$\overline{\Cal{P}M}_{0,k}$ near $\Sigma$.

We can see that there is an obvious partial order for the collection
of various topological types, $i.e.$ $I_1>I$ if the topological type
of $\Sigma_{I_1}$ can be obtained from $\Sigma_I$ by above gluing
procedure. Moreover, $\Sigma_{\alpha}\in \Cal{P}M^I_{0,k}$ if and
only if $t=0$ and $\tau=0$. Actually, we can get various curves in
$\Cal{P}M^{I_1}_{0,k}$ with $I_1>I$ by setting some of components of
$(t,\tau)$ be zero.

\bigskip

\noindent $\bullet$ {\it Compactness of stable moduli space }

\smallskip

We can define the Gromov-Floer topology, or for simplicity, say weak
topology on the moduli space $P{\Cal M}(\tilde{x}_-,\tilde{x}_+)$,
\begin{Definition}
A sequence $[\Sigma_n,v_n]=[(v_1,\cdots,v_K),(f_1,\cdots,f_l),o]$ or
$[V_n]$, $n=1,\cdots,\infty$, of stable connecting orbits is said to
Gromov-Floer weakly converge to a stable connecting orbit
$[\Sigma,v]$ $=[(v_1,\cdots,v_K)$, $(f_1,\cdots,f_l),o]_\infty$ or
$[V_\infty]$ if there exist representatives $V_n\in[V_n]$ and
$V_\infty\in[V_\infty]$ with domains $\Sigma_n$ and $\Sigma_\infty$
such that

$1^\circ$ $\Sigma_n\rightarrow\Sigma_\infty $, as
$n\rightarrow\infty$. This means that $\exists\ \Sigma'_n\in{\Cal
P}M_{0,k}^{I_{\Sigma_n}}$, $\Sigma'_\infty\in{\Cal
P}M_{0,k}^{I_{\Sigma_\infty}}$ with minimal marked points and
identification maps $\phi_n,\phi_\infty$ satisfying
$\Sigma_n=\phi_n(\Sigma'_n)$,
$\Sigma_\infty=\phi_\infty(\Sigma'_\infty)$ such that when $n$ is
sufficiently large $\Sigma'_n$ is in the neighborhood of
$\Sigma'_\infty$ and is represented by
$\Sigma'_{(\alpha_n,t_n,\tau_n)}$ with
$(\alpha_n,t_n,\tau_n)\rightarrow 0$.

$2^\circ$ For each compact set $K\subset
\Sigma'_\infty\setminus\{{\rm double\
points}\}\cup\{\cup_iL_{i,\pm\infty}\}$, if $n$ is sufficiently
large, let $K_n$ be the corresponding subset of $\Sigma'_n$ via
gluing in $\overline{\Cal{P}M}_{0,k}$, then $(V_n\circ\phi_n)|_K$ is
$C^\infty$-convergent to $(V_\infty\circ\phi_\infty)|_K$.

$3^\circ$ $\lim_{n\rightarrow\infty}E(V_n)=E(V_\infty)$.
\end{Definition}

Using the same method as in \cite{LT} (or in \cite{RT}
\cite{FO},{\it etc.}), we can prove

\begin{Theorem}
$P{\Cal M}(\tilde{x}_-,\tilde{x}_+)$ is Hausdorff and compact in the
sense of weak topology. Moreover, if
$\{[V_n]\}_{n=1}^{\infty}\rightarrow [V]_\infty$ in $P{\Cal
M}(\tilde{x}_-,\tilde{x}_+)$, then $E(V_n)\rightarrow E(V_\infty)$
and for sufficiently large $n$, $\mu(V_n,H)=\mu(V_\infty,H)$.
\end{Theorem}

\noindent {\it Remark}. Recall the Definition \ref{stable} we see
that each main component $v_i$ is the restriction of
$u_i=\pi(\tilde{u}_i)$, where $\tilde{u}_i\in {\Cal
M}(\tilde{x}_i,\tilde{x}_{i+1})$ is the lift of a
$(J,H)$-holomorphic map from the open domain $\R^2$, to the strip
$\R\times[0,1]$. The Gromov weak compactness arguments of Floer in
\cite{F3} are carried out for the special case $\phi=id$, and can be
easily generalized to arbitrary $\phi$ to prove the compactness of
the moduli space ${\Cal M}(\tilde{x}_-,\tilde{x}_+)$ modulo
splitting. Then for our moduli space of stable connecting orbits
$P{\Cal M}(\tilde{x}_-,\tilde{x}_+)$, we just need do the same
bubbling-off analysis as \cite{LT} or \cite{FO} restricting on the
subset $\R\times[0,1]$. Also we require that the bubble components
appear simultaneously on both sides of the boundary because of the
$\phi$-periodicity condition (\ref{002}). Here we will not repeat
the proof which can be found in the references listed above.

So the moduli space $P{\Cal M}(\tilde{x}_-,\tilde{x}_+)$ is the
stable compactification of the moduli space of connecting orbits
$\Cal{M}(\tilde{x}_-,\tilde{x}_+)$.

\subsection{Small neighborhood of stable moduli space}

Following \cite{LT}, with the difference in that our main components
are not cylinders but strips, we give a sketchy description of the
deformation of stable orbits under the topological change of their
domains. We firstly consider stable orbits with fixed intersection
pattern, then use the gluing procedure to deal with the stable
orbits with different intersection patterns. Then we get a
neighborhood $\Cal{W}$ of $P{\Cal M}$, which can be locally
uniformized and is a (partially smooth) orbifold.

Each stable orbit $(\Sigma,V)$ or $V$ consists of some main
components, denoted by $V_m$, and some bubble components, denoted by
$V_b$. Each $V$ determines an intersection pattern $D_V$ which is
determined by (I) the intersection pattern $I_{\Sigma}$ of the
domain $\Sigma$ and (II) the relative homotopy class of each main
component $V_m$ fixing its two end lines $L_{m,\pm\infty}$ and the
homology class represented by each bubble component $V_b$. Recall
the definition \ref{W^kp} we see that each
$V\in\Cal{B}(\tilde{x}_-,\tilde{x}_+)$ has the so-called {\it
effective} intersection pattern $D_V$ defined by Liu-Tian\cite{LT}.
In general, in the $W^{k,p}$-category, we say that an intersection
pattern $D$ is effective if $D=D_V$ with $V$ being a stable
$(J,H)$-map. So we can define the energy of each $V$ or $D$ by
$$E(V)=E(D=D_{V'})=E(V'),\ V'\in P\Cal{M},$$
where $E(V')$ is defined as (\ref{005}). If we denote the set of
intersection pattern with bounded energy by $D^e=\{D|E(D)\le e\}$,
then using the Lemma \ref{bound} we know that there are at most
finite number of marked points and ghost bubble components,
consequently, the set $D^e$ is finite.

Thus, according to intersection pattern, we can stratify
$P\Cal{M}=\cup_{D}P\Cal{M}^D$, where each
$$P\Cal{M}^D(\tilde{x}_-,\tilde{x}_+)=\{[V]|V\in
P\Cal{M}(\tilde{x}_-,\tilde{x}_+),\ D_V=D\},$$ and similarly, denote
the stratified set of $W^{k,p}$-orbits with bounded energy $E(V)\le
e$ by $\Cal{B}^e=\cup_{D}\Cal{B}^{e,D}$. Since for Given
$\tilde{x}_\pm$, the energy $E(V)$ of any stable connecting orbit
$[V]\in P\Cal{M}$ is bounded, if $e$ is sufficiently large, then
$P\Cal{M}(\tilde{x}_-,\tilde{x}_+)\subset
\Cal{B}^e(\tilde{x}_-,\tilde{x}_+)$. We will use $\Cal{B}^e$ as the
ambient space of $P\Cal{M}$, after taking an $e$ once for all, we
will omit the superscript and still denote the space by $\Cal{B}$.
We firstly study the space $\Cal{B}^D$ with fixed intersection
pattern $D$ and therefore $I_\Sigma$ of their domains.

Since the reparameterization group $G_\Sigma$ is noncompact, we have
no nice structure of $\Cal{B}^D$. While, we will show that near
$P\Cal{M}^D$ the action $G_\Sigma$ has a good slicing. To this end,
we choose a representative $V$ of $[V]\in P\Cal{M}^D$. Let $V_m,\
m=1,\cdots,M$ and $V_b,\ b=1,\cdots,B$ be its free main and bubble
components. For simplicity, we assume each main or bubble free
component has only one generic marked point, say $(0,\frac{1}{2})$
and $0$, respectively. We take locally a codimension 1 small disc
${\bf H}_m$ near $V_m(0,\frac{1}{2})$ for each free main component
$V_m$ such that $V_m|_{\R\times\{\frac{1}{2}\}}$ is transversal to
${\bf H}_m$, and a codimension 2 small disc ${\bf H}_b$ near
$V_b(0)$ for each free bubble component $V_b$ such that $V_b$ is
transversal to ${\bf H}_b$ at 0. Let $${\bf H}=\prod_{m=1}^M{\bf
H}_m\times\prod_{b=1}^B{\bf H}_b.$$

We define the distance in $\Cal{B}^D$ (each element has $K+L$
components and $d$ double points) as
$\|V-V'\|_{\Cal{B}^D}=\sum_{i=1}^{K+L}\|V_i-V'_i\|_{k,p}+\sum_{j=1}^ddist(z_j-z'_j)$,
where $dist$ is the distance function in the domain $\Sigma$ and the
Sobolev number $k,p$ will be taken carefully so that the $W^{k,p}$-
norm should be stronger than the $C^1$-topology. Now for a
sufficiently small $\epsilon$-neighborhood
$\widetilde{U}_{\epsilon}^D(V)=\{W|\|W-V\|_{\Cal{B}^D}\le\epsilon\}$
(where the $\epsilon$ is needed small enough so that
$[W(\Sigma_j)]=[V(\Sigma_j)]$ for all $j$), we can define a slicing
(at least with respect to those group actions of $\prod_{m=1}^PG_m$)
of $\widetilde{U}_{\epsilon}^D(V)$ as
$$\widetilde{U}_{\epsilon}^D(V,{\bf H})=\{W|W\in\widetilde{U}_{\epsilon}^D(V),\
W_m(0,\frac{1}{2})\in{\bf H}_m,\ W_b(0)\in{\bf H}_b\},$$ with taking
$m$ from $1$ to $M$ and $b$ from $1$ to $B$. So the problem is to
deal with the bubble component $V_b$.

Recall $G_\Sigma$ is the reparameterization group of stable maps. We
denote the automorphism group of $V$ by $\Gamma_V=\{g|g\in
G_\Sigma,\ V\circ g=V\}$. $\Gamma_V$ is a finite group since $V$ is
stable. And it is generated by the subgroup
$\prod_{m=1}^K\Gamma_{V_m}\times\prod_{b=1}^L\Gamma_{V_b}\times\Gamma_I$,
where $$\Gamma_{V_{m,b}}=\{g_{m,b}|V_{m,b}\circ g_{m,b}=V_{m,b},
{\rm preserving\ double\ points } \}$$ and $$\Gamma_I=\{g\in
G_\Sigma| {\rm \ interchanging\ components\ of\ }\Sigma {\rm \ and\
preserving}\ V \}.$$

If the automorphism group $\Gamma_V$ is trivial, for small enough
$\epsilon$, the projection $\pi_V:\widetilde{U}_{\epsilon}^D(V,{\bf
H})\rightarrow{\Cal W}^D$ which takes the point $W$ to its
equivalence class $[W]$ is injective, thus $[V]$ has a neighborhood
in ${\Cal W}^D$ modeled as an open subset in a Banach space. If
$\Gamma_V$ is nontrivial, that is to say, either some $V_b$ is a
multiple covering sphere or the automorphism interchanges components
of $\Sigma$, then as in section 2 of \cite{LT}, we can extend its
action to a linear action on $\widetilde{U}_{\epsilon}^D(V,{\bf H})$
in such a way that a neighborhood $U_{[V]}^D$ of $[V]$ in ${\Cal
W}^D$ can be identified with the quotient
$\widetilde{U}_{\epsilon}^D(V,{\bf H})/\Gamma_V$. For simplicity, we
only state how to extend the action of
$\widetilde{\Gamma}_V=\prod_{b=1}^L\Gamma_{V_b}=\Gamma_{V_b}$
($i.e.$ only one bubble component) to
$\widetilde{U}_{\epsilon}^D(V,{\bf H})$. The general case is in
principle same. We denote the set of pre-image points of the bubble
component $W_b$ by
$$W_b^{-1}(W_b(0))=\{y_1=0,y_2,\cdots,y_{n_b}\}.$$ Similar to the
Lemma 2.2 in \cite{LT} with the difference in that maybe some of
$y_i's$ is on the boundary of the domain, we can obtain the
analogous conclusion that for sufficiently small $\epsilon$ and
$\delta$, and for any $W\in \widetilde{U}_{\epsilon}^D(V,{\bf H})$,
there exist exactly $n_b$ points, $y_1(W_b),\cdots,y_{n_b}(W_b)$
such that for each $i$, $y_i(W_b)$ is in a $\delta$-disc or half
$\delta$-disc centering at $y_i$ (denoted by $D_{\delta}(y_i)$ or
$HD_{\delta}(y_i)$), and $$W_b^{-1}({\bf
H}_b)=\{y_1(W_b),\cdots,y_{n_b}(W_b)\}.$$

Let $g_i$ be the automorphism of $S^2$ such that $g_i(y_1)=y_i$,
$i=1,\cdots,n_b$, $g_i(1)=1$, $g_i(\infty)=\infty$, where we choose
$y_1=0$. For any $W\in \widetilde{U}_{\epsilon}^D(V,{\bf H})$, we
define an automorphism of $S^2$ as $$g_i^W: y_1\rightarrow y_i(W),\
1\rightarrow 1,\ \infty\rightarrow\infty.$$ Let
$r=\min_{i>m}\|V-V\circ g_i\|$. So if $\epsilon\ll\epsilon_1\ll r$,
then $W\circ g_i^W\in\widetilde{U}_{\epsilon_1}^D(V,{\bf H})$ if and
only if $i\le m$. This gives rise to an action of
$\widetilde{\Gamma}_V$ on $\widetilde{U}_{\epsilon_1}^D(V,{\bf H})$:
$$W\ast g=W\circ g^W,$$ for $W\in\widetilde{U}_{\epsilon_1}^D(V,{\bf
H})$, $g\in\widetilde{\Gamma}_V$.

We also see that for any give two elements $W_1$ and $W_2$ in
$W\in\widetilde{U}_{\epsilon_1}^D(V,{\bf H})$, $W_1$ and $W_2$ are
equivalent if and only if there exists a $g\in \Gamma_V$ such that
$W_1=W_2\ast g$. So if we replace
$\widetilde{U}_{\epsilon_1}^D(V,{\bf H})$ by the
$\Gamma_V$-invariant subset
$\cup_{g\in\Gamma_V}g(\widetilde{U}_{\epsilon}^D(V,{\bf H}))$, then
the action constructed above is a smooth right action on
$\cup_{g\in\Gamma_V}g(\widetilde{U}_{\epsilon}^D(V,{\bf H}))$, and a
neighborhood of $[V]$ in $\Cal{B}^D$ is homeomorphic to
$\cup_{g\in\Gamma_V}g(\widetilde{U}_{\epsilon}^D(V,{\bf
H}))/\Gamma_V$. For simplicity of notation, we still write the
$\Gamma_V$-invariant subset
$\cup_{g\in\Gamma_V}g(\widetilde{U}_{\epsilon}^D(V,{\bf H}))$ as
$\widetilde{U}_{\epsilon}^D(V,{\bf H})$ if no danger of confusion.

We say the triple $(\widetilde{U}_{\epsilon}^D(V,{\bf
H}),\Gamma_V,\pi_V)$ is a local uniformizer for $[V]$ in ${\Cal
W}^D$, where $\pi_V$ is the quotient map. In other words, there
exists a neighborhood ${\Cal W}^D$ of $P{\Cal M}^D$ in the space
${\Cal B}^D$ of all stable orbits with intersection pattern $D$
which is covered by local uniformizers
$(\widetilde{U}_{\epsilon}^D(V,{\bf H}),\Gamma_V,\pi_V)$. These
uniformizers give the neighborhood an orbifold structure.

Now the problem is to describe how the strata ${\Cal W}^D$ fit
together with topological change of the domain. We also need the
gluing procedure to describe the local information of a stable orbit
$V:\Sigma\rightarrow M$ in the full neighborhood ${\Cal W}$. Assume
the topological type of the domain $\Sigma$ is $I=I_{\Sigma}$.
Recall in the last subsection a $\delta$-neighborhood of the element
$\Sigma$ in ${\Cal P}M_{0,k}^{I_\Sigma}$ can be described by
parameters $\alpha=(\alpha_{l,k})$ with $\alpha_{l,k}\in
D_{\delta}(d_{l,k})$, and when the intersection pattern changes the
full $\delta$-neighborhood of $\Sigma$ is described by
$(\alpha,t,\tau)$ with $\|t\|,|\tau|\le\delta$, where $(t,\tau)$ are
the gluing parameters of the domain. Still from the delicate
pre-gluing procedure in \cite{LT} we can define the following stable
orbits as ``base point"
$$V_{\alpha}:\Sigma_{\alpha}\rightarrow M,\ \ \ \ \ \
V_{(\alpha,t,\tau)}:\Sigma_{(\alpha,t,\tau)}\rightarrow M$$ where
$V_{\alpha}$ is $W^{k,p}$-close to $V$ and $V_{(\alpha,t,\tau)}$ is
the pre-gluing of $V_{\alpha}$.

More precisely, we can define $V_{\alpha}=V\circ \psi_\alpha$, where
$\psi_\alpha : \Sigma_\alpha\rightarrow\Sigma$ is diffeomorphism
defined as follows. Fix a $r\ll\delta>0$, define $\psi_\alpha$ to be
identity on $\Sigma_\alpha\setminus \cup D_r(d_{l,k})$ and to be
rotation of $S^2$ on each $D_\delta(d_{l,k})$ bringing
$\alpha_{l,k}$ to $d_{l,k}$, where for simplicity $D_r$ and
$D_\delta$ denote half disc or disc according to whether $d_{l,k}$
is on the boundary or not. Since $r> 2\delta$, we can naturally
assume the image of $\psi_\alpha$ restricting to $D_\delta$ is
contained in $D_r$. So it is easy to smoothly extend $\psi_\alpha$
to all $\Sigma_\alpha$. And when $\delta$ is small enough,
$\psi_\alpha$ is smoothly close to identity. Thus, $V_\alpha$ is
$W^{k,p}$-close to $V$.

We can apply the very similar method in \cite{LT} to get the
pre-gluing $V_{(\alpha,t,\tau)}$ of $V_\alpha$. Recall that the
domain $\Sigma_{(\alpha,t,\tau)}$ can be derived from
$\Sigma_\alpha$ by gluing procedure listed in the last subsection.
The difference is in that we consider gluing strips $L_{i,+\infty}$
and $L_{i+1,-\infty}$ of each pair of connecting main components
instead of the annulus used by Liu-Tian. If there exist double point
$d_{l,k}$ in the boundary, there also will be no trouble, we still
can do the pre-gluing procedure firstly for a stable map defined on
a larger domain, then we can restrict the resulting pre-gluing map
to our original domain $\R\times[0,1]$.

Also we can easily show that there exists a map
$$\psi_{(\alpha,t,\tau)}: (\Sigma_{(\alpha,t,\tau)}, z^*_1,\cdots,
z^*_k)\rightarrow (\Sigma,z_1,\cdots,z_k),$$ which is injective
outside all small discs contain double points. Then the full
neighborhood of $V$ is denoted by $\widetilde{U}_{\epsilon}(V,{\bf
H})$, which contains all points $(\Sigma_{(\alpha,t,\tau)},
\hat{V},z'_1,\cdots,z'_k)$ satisfying that $\hat{V}$ is
$\epsilon$-close to $V\circ \psi_{(\alpha,t,\tau)}$ and each
$\psi_{(\alpha,t,\tau)}(z'_i)$ is $\epsilon$-close to $z_i$, where
the parameters $(\alpha,t,\tau)$ also vary in a
$\delta$-neighborhood.

Then as before we may extend the action of automorphism group
$\Gamma_V$ to the $\Gamma_V$-invariant set
$\cup_{g\in\Gamma_V}g(\widetilde{U}_{\epsilon}(V,{\bf H}))$, which
is still denoted by $\widetilde{U}_{\epsilon}(V,{\bf H})$, for
simplicity, if no danger of confusion. Let $U_{\epsilon}(V,{\bf
H})=\widetilde{U}_{\epsilon}(V,{\bf H})/\Gamma_V$, which can be
regarded as a small neighborhood in ${\Cal W}$. We just write ${\Cal
W}=\bigcup_{V\in P{M}}U_{\epsilon_V}(V,{\bf H})$. Moreover, the
$W^{k,p}$-topology of ${\Cal W}$ can be generated by
$U_{\epsilon}(V,{\bf H})$. This gives the ${\Cal W}$ an orbifold
structure (c.f. Lemma 2.6 in \cite{LT}). Also Liu-Tian proved that
the so-defined $W^{k,p}$-topology is equivalent to the Floer-Gromov
weak topology. This implies $P{\Cal M}$ is also Hausdorff and
compact with respect to the strong $W^{k,p}$ topology.

Consequently, we can take a finite union of the covering of $\Cal W$
as $$\{U_i=U_{\epsilon_{V_i}}(V_i,{\bf H}), i=1,\cdots,w \},$$ and
we use $\tilde{U}_i$ to denote its uniformizer with covering group
$\Gamma_i$.

Then we can locally define orbifold bundles ${\Cal L}_i$ over
$\tilde{U}_i$. For each $[V]\in U_i$, the fiber ${\Cal L}_i|_{[V]}$
over $[V]$ consists of all elements of
$L^{k-1,p}(\Lambda^{0,1}(V^*TM)),\ V\in[V]$ modulo equivalence
relation induced by pull-back of sections coming from identificaion
of the domains of $V_i$, where $\Lambda^{0,1}(V^*TM))$ is the bundle
of $(0,1)$-forms on $\Sigma$ with respect to the complex structure
on $\Sigma$ and the given compatible almost complex structure $J$ on
$(M,\omega)$. Then the local uniformizer $\tilde{\Cal L}_i$ of
$\Cal{L}_i$ is given by the union of
$L^{k-1,p}(\Lambda^{0,1}(\tilde{V}_i^*TM))$, $\tilde{V}_i\in
\tilde{U}_i$. The $\Gamma_i$ also acts on $\tilde{\Cal L}_i$ so that
${\Cal L}_i=\tilde{\Cal L}_i/\Gamma_i$. In this way, we can
reinterpret the $\bar{\partial}_{J,H}$-operator as a collection of
$\Gamma_i$-equivariant sections of these local orbifold bundles
$({\Cal L}_i,U_i)$.

More precisely, we will describe the construction in the rest of the
section. For each $W\in \widetilde{U}_\epsilon^D(V,{\bf H})$ or
$\widetilde{U}_\epsilon(V,{\bf H})$, the fiber is
$$\widetilde{\Cal L}^D(V)|_W=\widetilde{\Cal L}(V)|_W=\{\xi|\xi\in L^{k-1,p}(\Lambda^{0,1}(W^*TM))\},$$
where the $L^{k-1,p}$-norm is measured with respect to the metric on
the domain $\Sigma_W=\Sigma_{(\alpha,t,\tau)}$ induced by the gluing
construction from the metric on $\Sigma$ that is ``spherical" on
$\Sigma_b$ and flat on $\Sigma_m$.

For fixed intersection pattern $D$, $\widetilde{\Cal L}^D(V)$ is a
locally trivial Banach bundle over $\widetilde{U}_\epsilon^D(V,{\bf
H})$, and $\widetilde{\Cal L}(V)$ is locally trivial only when
restricted to each stratum $\widetilde{U}_\epsilon^D(V,{\bf H})$ of
$\widetilde{U}_\epsilon(V,{\bf H})$. So the topology of the bundle
$\widetilde{\Cal L}(V)$, when restricted to each stratum of
$\widetilde{U}_\epsilon(V,{\bf H})$, is well-defined. We will not
specify the topology right now, that will be done in the later
gluing construction of virtual cycle. Instead, now we consider the
more relevant following ``sub-bundle" $\widetilde{\Cal
L}^D_\delta(V)$ of $\widetilde{\Cal L}^D(V)$ without singularity,
defined as $$\widetilde{\Cal
L}^D_\delta(V)|_W=\{\xi|\xi\in\widetilde{\Cal L}^D(V)|_W,\ \xi=0\
{\rm on\ each}\ D_\delta(d_{l,k}) .\}$$

Let $\widetilde{U}_{\epsilon,\delta}^D(V,{\bf
H})=\bar{\partial}^{-1}_{J,H}(\widetilde{\Cal L}^D_\delta(V))$, then
we get a restricted bundle $$\widetilde{\Cal
L}^D_\delta(V)\rightarrow \widetilde{U}_{\epsilon,\delta}^D(V,{\bf
H}).$$ If $D\le D_1$, then we can move the fiber of $\widetilde{\Cal
L}^D_\delta(V)$ over some point in
$\widetilde{U}_{\epsilon,\delta}^D(V,{\bf H})$, by parallel
transformation, into the fiber of $\widetilde{\Cal
L}^{D_1}_{\delta_1}(V)$ over a neighborhood of the given point in
$\widetilde{U}_{\epsilon,\delta_1}^{D_1}(V,{\bf H})$, when
$\delta_1\ll\delta$. All these parallel transformations induce a
topology for the union $$\widetilde{\Cal
L}^0(V)=\cup_{D,\delta}\widetilde{\Cal
L}^D_\delta(V)\rightarrow\widetilde{U}_{\epsilon}^0(V,{\bf
H})=\cup_{D,\delta}\widetilde{U}_{\epsilon,\delta}^D(V,{\bf H}).$$

The $\Gamma_V$-actions on $\widetilde{U}_{\epsilon}^D(V,{\bf H})$
and $\widetilde{U}_\epsilon(V,{\bf H})$ can be lifted to the bundles
via pull-back. Recall $U^D_\epsilon(V,{\bf
H})=\widetilde{U}^D_\epsilon(V,{\bf H})/\Gamma_V$ and
$U_\epsilon(V,{\bf H})=\widetilde{U}_\epsilon(V,{\bf H})/\Gamma_V$,
we denote the orbifold bundles over them by ${\Cal
L}^D(V)=\widetilde{\Cal L}^D(V)/\Gamma_V$ and ${\Cal
L}(V)=\widetilde{\Cal L}(V)/\Gamma_V$, respectively.

With the above construction can see that the
$\bar{\partial}_{J,H}$-operator gives rise to a
$\Gamma_V$-equivariant section of the bundle $\widetilde{\Cal
L}(V)\rightarrow \widetilde{U}_\epsilon(V,{\bf H})$, it is smooth on
each stratum $\widetilde{U}^D_\epsilon(V,{\bf H})$, and continuous
when restricted to $\widetilde{\Cal L}^0(V)\rightarrow
\widetilde{U}_\epsilon^0(V,{\bf H})$. We still denote the section by
$\bar{\partial}_{J,H}$. The zero sets $\bar{\partial}_{J,H}^{-1}(0)$
in $\widetilde{U}^D_\epsilon(V,{\bf H})$ and
$\widetilde{U}_\epsilon(V,{\bf H})$ are just
$$\widetilde{U}^D_\epsilon(V,{\bf H})\cap P{\Cal
M}^D(J,H,\tilde{x}_-,\tilde{x}_+)$$ and
$$\widetilde{U}_\epsilon(V,{\bf H})\cap P{\Cal
M}(J,H,\tilde{x}_-,\tilde{x}_+).$$

Thus, we consider ${\Cal W}$ as a {\it multi-fold} which is a
partially smooth space---space with two topologies, which will be
defined in the next section.

\section{Abstract settings}

All arguments in this section can be found in \cite{LT} and
\cite{M99}, just for reader's convenience we give some related
definitions and notations used later .

\subsection{Partially smooth space and branched pseudomanifold}

\begin{Definition}
A Hausdorff space $Y$ is said to be partially smooth if it is the
image of a continuous bijection $i_Y:Y_{sm}\rightarrow Y$, where
$Y_{sm}$ is a finite union of open disjoint subsets, each of which
is a smooth Banach manifold.
\end{Definition}
We consider the collection of all partially smooth spaces as objects
of a category, and a morphism is a continuous map $f:Y\rightarrow X$
between two objects such that the induced map $f:Y_{sm}\rightarrow
X_{sm}$ is smooth, say a partially smooth map. We see such $Y$ is
stratified, the strata that are open subsets of $Y$ are called top
strata.

\begin{Definition}
A pseudomanifold of dimension $d$ is a compact partially smooth
space $Y$ such that a component $Y_{sm}$ is an oriented smooth
$d$-dimensional manifold which is mapped by $i_Y$ onto a dense open
subset $Y^{top}$ of $Y$, and the dimensions of all other components
of $Y_{sm}$ are no larger than $d-2$.
\end{Definition}
We denote by $Y^{sing}=Y-Y^{top}$ for the image of those lower
dimensional submanifolds. The following object is more general.

\begin{Definition}
A branched pseudomanifold $Y$ of dimension $d$ (without boundary) is
a compact partially smooth space such that its components have at
most dimension $d$. We denote the components of dimension $d$ by
$M_i$, those of dimension $d-1$ by $B_j$, and write
$$Y^{top}=\bigcup_iM_i,\ \ \ B=\bigcup_jB_j,\ \ \
Y^{sing}=Y-(Y^{top}\cup B).$$ Especially, we assume that for each
$i$ the set $M_i\cup_{j\in J_i}B_j$ ($j\in J_i$ if the closure of
$M_i$ in $Y$ meets $B_j$) has the structure of smooth manifold with
boundary which is compatible with its two topologies. And we assume
$\bar{B}_j-B_j\subset Y^{sing}$. We call $B$ the branched locus. A
branched $d$-pseudomanifold with boundary is defined similarly, but
some $d-1$ dimensional components of $\partial M_i$ are not any of
the branched locus. We denote the union of such $d-1$ dimensional
components by $\partial Y$.
\end{Definition}

In order to construct the virtual cycle, we need a suitable labeling
of the top components.
\begin{Definition}\label{label}
We say a branched pseudomanifold $Y$ (with or without boundary) is
labeled if its top components $M_i$ are oriented and have positive
rational labeling $\lambda_i\in \Q $ satisfying that for each $x\in
B$, if we pick an orientation of $T_xB$, and divide the components
$M_i$ which have $x$ in their closure into two groups $I^+,I^-$
according to whether the chosen orientation on $T_xB$ agrees with
the boundary orientation, then $\sum_{i\in I^+}\lambda_i=\sum_{j\in
I^-}\lambda_j$.
\end{Definition}
In particular, if $Y$ is of dimension zero, it is a collection of
oriented labeled points, and there is no compatibility condition at
the branch locus since it is empty.

If a compact branched pseudomanifold $Y$ of dimension $d$ is labeled
as above, then McDuff \cite{M99} in the following Lemma showed that
it can be regarded as a (relative) cycle which represents a
(relative) rational class.
\begin{Lemma}\label{5.1}
Let $Y$ be a (closed, oriented) branched and labelled pseudomanifold
of dimension $d$. Then every partially smooth map $f$ from $Y$ to a
closed manifold $X$ defines a rational class $f_*([Y])\in H_d(X)$.
\end{Lemma}
Proof. Let $Z$ be a smooth manifold of dimension complementary to
$Y$ and $g:Z\rightarrow X$ be a smooth map. Then we can jiggle $g$
so that it doesn't meet $f(Y^{sing})$ and so that it meets $\cup
f(M_i)$ transversally in a finite number of points. We then define
$$f\cdot g=\sum_i\delta_i\lambda_i,$$ where $\delta_i=+1$ or $-1$
for $i\in I^+$ or $I^-$. One can check this number is independent of
jiggling. Then we can define $f_*([Y])$ to be the unique rational
class such that the intersection number $$f_*([Y])\cdot
g_*([Z])=f\cdot g,\ \ \ {\rm for\ all}\ g:Z\rightarrow X.$$ So we
can look $Y$ as a cycle. \qed

If $Y$ is of dimension 1 with boundary, then from the way of
labeling we see that the oriented number of its boundary is zero.
More precisely, let $x\in\partial Y\cap M_i$, and denote by
$\overrightarrow{v}_i\in T_xM$ the outward unit normal vector. This
vector together with the orientation of $M_i$ determines a sign
$$\delta_i=\left\{
\begin{array}{ccc}
+1,\ \ {\rm if}\ \overrightarrow{v}_i\ {\rm is\ positively\ oriented}, \\

-1,\ \ {\rm if}\ \overrightarrow{v}_i\ {\rm is\ negatively\
oriented}.
\end{array}
\right.$$

Then for each boundary point we define a number
$$\rho(x)=\sum_{x\in M_i}\delta_i\lambda_i.$$
D. Salamon proved the following result \cite{S}
\begin{Lemma}\label{5.2}
Let $(Y,\lambda)$ be a compact oriented, branched, and labeled
1-pseudomanifold with boundary. For each $x\in\partial M$ we have a
rational number $\rho(x)$ defined above. Then $$\sum_{x\in\partial
M}\rho(x)=0.$$
\end{Lemma}

Another point is that, as much as possible, we are trying to avoid
specifying exactly how the strata of $Y_{sm}$ fit together. If $Y$
is branched then one does need some information of the normal
structure to the codimension 1 components, but one can often get
away without any other restrictions. However, later we will need to
consider the intersection of two partially smooth spaces, and in
order for this to be well-behaved one does need more structure. This
is the reason for the following definition.
\begin{Definition}\label{cone}
We will say that a partially smooth space $Y$ has normal cones if
every stratum $S$ in $Y$ has a neighborhood $N(S)$ in $Y$ whose
induced stratification is that of a cone bundle over a link. More
precisely, there is a commutative diagram
\begin{eqnarray}
\ \ N(S)_{sm}       & \longrightarrow & N(S)\nonumber \\
\pi'\downarrow  &             &\pi\downarrow  \\
\ \ S_{sm}     &\longrightarrow &\ \   S  \nonumber
\end{eqnarray}
where the maps $\pi'$ and $\pi$ are oriented locally trivial
fibrations with fiber equal to a cone over a link $L$. Here $L$ is
partially smooth, and the cone $C(L)$ is just the quotient $L\times
[0, 1]/L\times {0}$, stratifed so that it is the union of the vertex
(the image of $L\times {0}$) with the product strata in
$L_{sm}\times (0, 1]$. Moreover each stratum in $N(S)_{sm}$ projects
onto the whole of S and so is a locally trivial bundle over $S$ with
fiber equal to a component of $C(S)_{sm}$. In particular, we
identify S with the section of $N(S)$ given by the vertices of the
cones.
\end{Definition}

This definition implies that any stratum $S'$ whose closure
intersects $S$ must contain the whole of $S$ in its closure. In
fact, near $S$ it must look like the cone over some stratum of
$L_{sm}$. Any finite dimensional Whitney stratified space has this
normal structure. For example, any finite dimensional orbifold has a
stratification such that the isomorphism class of the automorphism
group $\Gamma_x$ at the point $x$ is constant as $x$ varies over
each stratum, and it is easy to see that this stratification has
normal cones as defined here. It is also not hard to use arguments
similar to those in last subsection to show that the neighborhood
${\Cal W}$ of $P{\Cal M}$ in the space of all stable maps has normal
cones with respect to the fine stratification.

\subsection{Multi-fold, multi-bundle, and multi-section}

Now still following \cite{M99}, we introduce the concept of
multi-fold which is a generalization of orbifold. It can be
considered as atlas (or covering) of a space $\Cal{W}$ with two
topologies that locally is an orbifold in the partially smooth
category, in which the inclusions that relate one uniformizer to
another in an orbifold are replaced by {\it fiber products} used by
Liu-Tian, also they consider the fiber product as a certain kind of
``global uniformizer". We then show the definition of multi-sections
of a multi-bundle over a multi-fold. All maps, spaces and group
actions considered below are in this partially smooth category.

Suppose that a space $\Cal{W}$  is a finite union of open sets
$U_i$, $i = 1,\cdots  w$, each with uniformizers
$(\tilde{U}_i,\Gamma_i,\pi_i )$ with the following properties. Each
$\Gamma_i$ is a finite group acting on $\tilde{U}_i$ and the
projection $\pi_i$ is the composite of the quotient map
$\tilde{U}_i\rightarrow \tilde{U}_i/\Gamma_i$ with an identification
$\tilde{U}_i/\Gamma_i=U_i$. The inverse image in $\tilde{U}_i$ of
each stratum in $U_i$ is an open subset of a (complex) Banach space
on which $\Gamma_i$ acts complex linearly. For simplicity, we will
assume that $\Gamma_i$ acts freely on the points of the top strata
in $\tilde{U}_i$, and that the isomorphism class of the stabilizer
subgroup Stab$_i(\tilde{x}_i)$ of $\Gamma_i$ is fixed as
$\tilde{x}_i$ varies over a stratum of $\tilde{U}_i$, where the
Stab$_i(\tilde{x}_i)$ is the subgroup of $\Gamma_i$ that fixes
$\tilde{x}_i$.

For each subset $I=\{i_1,\cdots,i_p\}$ of $\{1,\cdots,w\}$ set
$U_I=\cap_{j\in I}U_j$, and $U_{\emptyset}=\emptyset$. Let ${\Cal
N}=\{I|\ U_I\neq\emptyset\}$. For $I\in {\Cal N}$, denote
$\Gamma_I=\prod_{j\in I}\Gamma_j$.
\begin{Definition}\label{product}
The fiber product of those $\tilde{U}_j$, $j\in I$ is
$$\tilde{U}_I=\{ \tilde{x}_I=(\tilde{x}_j)_{j\in I}|\ \pi_j(\tilde{x}_j)=\pi_l(\tilde{x}_l)\in U_I, {\rm for\
all}\ j,l\in I \}$$ which is contained in $\prod_{j\in
I}\tilde{U}_j$ with the two topologies induced from $\prod_{j\in
I}\tilde{U}_j$.
\end{Definition}
Since it is easy to see $\Gamma_I$ acts on $\tilde{U}_I$, the
quotient $\tilde{U}_I/\Gamma_I\simeq U_I$. Denote the projection by
$\pi_I:\tilde{U}_I\rightarrow U_I$. Note that $\tilde{U}_I,\ U_I$
are also  partially smooth spaces. And the isomorphism class of the
stabilizer subgroup
$${\rm Stab}_I(\tilde{x}_I)=\prod_{j\in
I}{\rm Stab}_j(\tilde{x}_j)$$
of $\tilde{x}_I$ in $\Gamma_I$ is
constant on each stratum, and trivial at points of top strata.
Roughly speaking, we can consider the fiber product as a substitute
for $\bigcap_{j\in I}\tilde{U}_j$. The topological structure can be
understood in terms of notions {\it local component} and {\it
desingularization} defined by Liu-Tian \cite{LT2}.

More precisely, Given a point $\tilde{x}_I\in\tilde{U}_I$, choose
$i_0\in I$ and a small open neighborhood $\tilde{N}$ of
$\tilde{x}_{i_0}\in\tilde{U}_{i_0}$. Then a neighborhood of
$\tilde{x}_I$ in $\tilde{U}_I$ can be identified with the set
$$\{(\Upsilon_j\circ\iota_{ji_0}(\tilde{y}))_{i\in I}:\ \tilde{y}\in\tilde{N},\ \Upsilon_{i_0}=id,
\ \Upsilon_j\in{\rm Stab}_j(\iota_{ji_0}(\tilde{y})),\ j\neq
i_0\}.$$ It is a finite union of sets
$\tilde{N}_\Upsilon(\tilde{x}_I)$, where each element $\Upsilon$ is
in the group $${\rm Stab}_{I-i_0}(\tilde{x}_I)=\prod_{j\in
I-i_0}{\rm Stab}_j(\tilde{x}_j).$$ In the partially smooth category,
each such set is homeomorphic to $\tilde{N}$. It is clear that the
germs at $\tilde{x}_I$ of the sets $\tilde{N}_\Upsilon(\tilde{x}_I)$
are independent of the choices of $i_0$ and $\tilde{N}$, as is the
isomorphism class of the ``reduced" group ${\rm
Stab}'_I(\tilde{x}_I)={\rm Stab}_{I-i_0}(\tilde{x}_I)$. These germs
are called the local components of $\tilde{U}_I$ at $\tilde{x}_I$,
denoted by $\langle\tilde{N}_\Upsilon(\tilde{x}_I)\rangle$. Since
the points $\tilde{x}_I$ in the top strata have trivial stabilizer
groups, they only have a single local component. Then the
``desingularization" $\widehat{U}_I$ of $\tilde{U}_I$ is defined to
be the union of all such local components
$$\widehat{U}_I=\{(\tilde{x}_I,\langle\tilde{N}_\Upsilon(\tilde{x}_I)\rangle):\ \tilde{x}_I\in\tilde{U}_I,
\ \Upsilon\in{\rm Stab}'_I(\tilde{x}_I) \}.$$ We can topologize the
desingularization $\widehat{U}_I$ so that a germ of neighborhood
contained in $\widehat{U}_I$ at the point
$(\tilde{x}_I,\langle\tilde{N}_\Upsilon(\tilde{x}_I)\rangle)$ is
homeomorphic to the local component
$\langle\tilde{N}_\Upsilon(\tilde{x}_I)\rangle$ itself. So the
projection $$proj:\ \widehat{U}_I\rightarrow\tilde{U}_I$$ is locally
a homeomorphism onto its image. Since locally $\widehat{U}_I$ is
homeomorphic to the initial sets $\tilde{U}_i$, the extra
singularities of $\tilde{U}_I$, introduced by constructing the fiber
product, are of no trouble.

If $J=(j_1,\cdots,j_q)\subset I=(i_1,\cdots,i_p)$, there are two
natural projections
$$\pi_J^I:\tilde{U}_I\rightarrow \tilde{U}_J=\tilde{U}_I/\Gamma_{I-J},\ \ \lambda_J^I:\Gamma_I\rightarrow\Gamma_J,$$
induced from the corresponding projection $\prod_{i_k\in
I}\tilde{U}_{i_k}\rightarrow\prod_{j_l\in J}\tilde{U}_{j_l}$ such
that $\pi_J\circ\pi_J^I=\iota_J^I\circ\pi_I$, where $\iota_J^I$ is
the inclusion $U_I\hookrightarrow U_J$. We see that if
$\pi_J^I(\tilde{y}_I)=\tilde{x}_J$ is in a top stratum, then
$(\pi_J^I)^{-1}$ has $|\Gamma_{I-J}|$ points.

Let $\{{\bf V}_I\}$ be a open cover of ${\Cal W}$ such that ${\bf
V}_I\subset U_I$ for each $I$. From the example 4.10 in \cite{M99}
we know that in general the sets ${\bf V}_j,\ j=1,\cdots,w,$ no
longer cover $\Cal{W}$. Then we define $\tilde{{\bf
V}}_I\subset\tilde{U}_I$ to be the inverse image of ${\bf V}_I$
under the map $\pi:\tilde{U}_I\rightarrow U_I$. So $\Gamma_I$ acts
on $\tilde{{\bf V}}_I$ and we still write the quotient map as
$\pi:\tilde{{\bf V}}_I\rightarrow {\bf V}_I\simeq\tilde{{\bf
V}}_I/\Gamma_I$. We define the projection $\pi_J^I:\tilde{{\bf
V}}_I\rightarrow \tilde{\bf V}_J$ as the restriction of the
projection $\pi_J^I$ above with domain $(\pi_J^I)^{-1}(\tilde{{\bf
V}}_J)$.

\begin{Definition}\label{m-fold}
A multi-fold atlas for $\Cal{W}$ is a collection $$\widetilde{\Cal
V}=\{(\tilde{{\bf V}}_I,\Gamma_I,\pi_J^I,\lambda_J^I),\ I\in {\Cal
N}\}.$$ The $\Cal{W}$ with such an atlas is called a multi-fold.
\end{Definition}
The motivation of considering such an atlas $\widetilde{\Cal V}$ or
a subcover $\{{\bf V}_I\}$ of $\{U_I\}$ is in that in general
(especially, in our application that ${\Cal W}$ is a neighborhood of
stable moduli space $P{\Cal M}$) the suitable chosen $\{{\bf V}_i\}$
will have simpler overlaps rather than $\{U_i\}$ and when the sets
$U_i$ overlap too much there are no non-equivariant multi-sections
which will be defined below. Since perturbation in the class of
equivariant sections is not sufficient to realize regularity, we
hope to obtain non-equivariant ones.

We now introduce the concept of {\it multi-bundle}. Let
$\Cal{E}=\bigcup_iE_i $ be a space with two topologies, and a map
$p:\Cal{E}\rightarrow \Cal{W}$ with the property that each set
$E_i=p^{-1}(U_i)$ has a local uniformizer
$(\tilde{E}_i,\Gamma_i,\Pi_i)$ such that the following diagram
commutes
\begin{eqnarray}
\ \ \tilde{E}_i       & \overrightarrow{\Pi_i} & E_i \nonumber \\
\tilde{p} \downarrow  &             & p \downarrow  \\
\ \ \tilde{U}_i     &\overrightarrow{\pi_i} &\ \   U_i.  \nonumber
\end{eqnarray}
And the map $\tilde{p}:\tilde{E}_i\rightarrow\tilde{U}_i$ is
required to be $\Gamma_i$-equivariant such that its restriction to
each stratum of $\tilde{U}_i$ is a locally trivial vector bundle.
Then the fiber $\tilde{F}(\tilde{x}_i)$ of $\tilde{p}$ at each point
$\tilde{x}_i$ is a vector space, but with no natural identification
of two different fibers if they are over points in different strata.
We denote by $\tilde{E}_I$ the restriction to $\tilde{\bf V}_I$ of
the fiber product of the $\tilde{E}_i,\ i\in I$, over $\tilde{U}_I$.

We suppose that the orbifold structure on ${\Cal W}$ can lift to one
on ${\Cal E}$, then ${\Cal E}$ had the same local structure as
${\Cal W}$. As above, we define a multi-fold atlas $\widetilde{\Cal
E}=\{\tilde{E}_I\}$ for ${\Cal E}$. More precisely, the elements of
the fiber $\tilde{F}(\tilde{x}_I)$ of $\tilde{E}_I$ at
$\tilde{x}_I\in\tilde{V}_I$ are $(\tilde{x}_i,\tilde{v}_i)_{i\in
I}$, where $\tilde{v}_i\in\tilde{F}(\tilde{x}_i)$, and for all
$j,l\in I$, $$\Pi_j(\tilde{v}_j)=\Pi_l(\tilde{v}_l)\in E_I.$$ Also
we can get the desingularization
$$\widehat{E}_I\rightarrow\widehat{\bf V}_I$$ of
$\tilde{E}_I\rightarrow\tilde{\bf V}_I$, which is an honest vector
bundle since locally it has the same structure as the maps
$\tilde{E}_i\rightarrow\tilde{U}_i$, $i\in I$. Thus
$\tilde{E}_I\rightarrow \tilde{V}_I$ is a finite union of vector
bundles with each fiber being a finite union of vector spaces.
\begin{Definition}
We say that the map $\tilde{p}:\ \widetilde{\Cal
E}\rightarrow\widetilde{\Cal V}$ constructed above is a multi-bundle
with respect to the multi-fold atlas defined in the Definition
\ref{m-fold}.
\end{Definition}

Now, we come to the definition of {\bf  multi-section} of the
multi-bundle $\tilde{p}:\ \widetilde{\Cal
E}\rightarrow\widetilde{\Cal V}$. Intuitively we consider it as a
{\it compatible collection} $\{\tilde{s}_I\}$ of multi-valued
sections. That is to say, we have a nonempty finite subset
$\tilde{s}_I(\tilde{x}_I)$ of $\tilde{F}_I(\tilde{x}_I)$ for each
point $\tilde{x}_I\in\tilde{{\bf V}}_I$. In our application, it is
enough to consider multi-sections that are single-valued when lifted
to the desingularization $\widehat{E}_I\rightarrow\widehat{\bf
V}_I$. We assume that there exist sections $\hat{s}_I:\widehat{\bf
V}_I\rightarrow\widehat{E}_I$ such that
$$\tilde{s}_I(\tilde{x}_I)=\{proj\circ\hat{s}_I(\hat{x}_I,\langle \tilde{N}_\Upsilon(\hat{x}_I)\rangle):
\ \Upsilon\in {\rm Stab}'_I(\hat{x}_I) \}.$$ Note that when
$I=\{j\}$ the section $\tilde{s}_j=\hat{s}_j$ is a single-valued
section of the bundle $\tilde{E}_j\rightarrow\tilde{U}_j$ which is
maybe non-equivariant. And when $\tilde{x}_I$ is in the top stratum
the set $\tilde{s}_I(\tilde{x}_I)$ has only one element in
$\tilde{F}_I(\tilde{x}_I)$.

Let $J\subset I$, recall there is projection $\Pi_J^I:\
\tilde{E}_I\rightarrow\tilde{E}_J$. For a multi-section
$\tilde{s}_J$ we define its pullback $(\Pi_J^I)^{-1}(\tilde{s}_J)$
to be a multi-valued section so that each $\tilde{x}_I$ is
associated with the full inverse image
$(\Pi_J^I)^{-1}(\tilde{s}_J(\tilde{x}_J))$ where
$\tilde{x}_J=\pi_J^I(\tilde{x}_I)$. Then the {\it compatibility}
condition requires that
\begin{equation}\label{compatibility}
\tilde{s}_I|_{(\pi_J^I)^{-1}(\tilde{{\bf
V}}_J)}=(\Pi_J^I)^*\tilde{s}_J.
\end{equation}

It is easy to see that if for all local bundles $\tilde{p}:\
\tilde{E}_i\rightarrow\tilde{U}_i$ a family of
$\Gamma_i$-equivariant sections $\tilde{s}_i$ are compatible, then
the multi-bundle $\widetilde{\Cal E}\rightarrow\widetilde{\Cal V}$
has a multi-section. In our application the $\bar{\partial}_{J,H}$
operator is such an example. However, because of the appearance of
multi-covered spheres with negative Chern class, perturbations of
the Cauchy-Riemann equations in the equivariant class can not
regularize the (stable) moduli space. Therefore, we have to consider
a method of extending non-equivariant multi-valued sections of
bundles $\tilde{E}_i\rightarrow\tilde{U}_i$ to multi-sections of
$\widetilde{\Cal E}\rightarrow\widetilde{\Cal V}$, such that they
can be regarded as global perturbations of the
$\bar{\partial}_{J,H}$ operator. Actually, if we choose a ``good"
subcover $\{{\bf V}_I\}$ of $\{U_I\}$ which do not overlap too much,
the extension is possible. The following lemma shows the existence
of such a good $\{{\bf V}_I\}$, it is taken from \cite{LT} and
\cite{M99}. The notation $A\subset\subset B$ means that the closure
of $A$ is contained in $B$.
\begin{Lemma}\label{VI}
Given a finite open covering $\{U_i,i=1,\cdots,N\}$ of a compact
subset $P{\Cal M}$ of $\Cal W$ and $U_I$ is defined as above, then
there are open subsets $U^0_i\subset\subset U_i$ and ${\bf
V}_I\subset U_I$ satisfying (i)$P{\Cal M}\subset\cup_iU_i^0$,
$P{\Cal M}\subset\cup_I{\bf V}_I$; (ii) If $i\notin \ I$, $U_i^0\cap
{\bf V}_I=\emptyset$; (iii) if ${\bf V}_I\cap {\bf
V}_J\neq\emptyset$ then $I\subset J$ or $J\subset I$.
\end{Lemma}
Proof. For $n=0,1,\cdots,N$ choose open coverings $\{U^n_i\}$,
$\{W^n_i\}$ of $P{\Cal M}$ satisfying
$$U_i^0\subset\subset W_i^1 \subset\subset U_i^1 \subset\subset W_i^2 \subset\subset\cdots\subset\subset U_i^N=U_i.$$
Then, let $\kappa=|I|$, we define $${\bf
V}_I=W^{\kappa}_I-\bigcup_{J:|J|>\kappa}{\rm Closure\ of}\
U^{\kappa+1}_J.$$ All properties hold obviously.\qed

By shrinking $\Cal W$, we may suppose ${\Cal W}=\cup_i U_i^0=\cup_I
{\bf V}_I$. For a multi-bundle $\widetilde{\Cal
E}\rightarrow\widetilde{\Cal V}$ as above, we assume that for some
$j$ we have a section $\sigma(j):\tilde{U}_j\rightarrow\tilde{E}_j$
of the bundle $\tilde{E}_j\rightarrow\tilde{U}_j$ with support in
$\tilde{U}^0_j$. Then, for each $I$ if $j\notin\ I$, we can define
$\tilde{s}(j)_I$ to be the zero section of
$\tilde{E}_I\rightarrow\tilde{{\bf V}}_I$ and, otherwise, we can
define $\tilde{s}(j)_I$ to be the restriction to $\tilde{{\bf V}}_I$
of the pullback to $\tilde{U}_I$ of the graph of $\sigma(j)$. This
induces a multi-section $\tilde{s}(j)$ of $\widetilde{\Cal
E}\rightarrow\widetilde{\Cal V}$, $j=1,2,\cdots,w$. That is to say,
choosing suitable covers $\{U_i^0\}$ and $\{{\bf V}_I\}$ as in the
lemma \ref{VI}, it is always possible to construct a multi-section.
In applications, we need the section $\sigma(j)$ satisfies some
generic conditions and require the boundary of the support of
$\sigma(j)$ is a union of strata (c.f. \cite{M99}). In the next
section we will use a form of
$\tilde{s}=\bar{\partial}_{J,H}+\sum_{j}\tilde{s}(j)$ with generic
perturbation term to define a Fredholm multi-section which is
transversal to zero section and derive the virtual moduli cycle.

\subsection{Constructing the virtual cycle}\label{construct}
Here, we will  construct  a branched pseudomanifold from the
multi-section defined above, that is just the virtual cycle.

In the following, we always take the cover $\{{\bf V}_I\}$ as in the
lemma \ref{VI}. In the partially smooth category we also have a
similar definition of {\it Fredholm system} as the Definition 1.1.
For each $I$, we denote by $gr(\tilde{s}_I)$ for the graph of the
section which is the union $\bigcup_{\tilde{x}_I\in \tilde{{\bf
V}}_I}\tilde{s}_I(\tilde{x}_I)$, we always require the graph is an
object and the projection $gr(\tilde{s}_I)\rightarrow\tilde{{\bf
V}}_I$ is a morphism in the partially smooth category. This amounts
to requiring that for each $j\in I$ the $gr(\tilde{s}_j)$ has a
stratification which is compatible with the projection to
$\tilde{{\bf V}}_j$ and is preserved by the action of the group
$\Gamma_j$ such that the rank of the stabilizer is constant on
strata (after refining of the stratification of $\tilde{{\bf
V}}_I$). Let $\widetilde{\Cal E},\widetilde{\Cal V},\Cal{W}$ are
defined as above, $\tilde{s}$ is a multi-section of $\widetilde{\Cal
E}\rightarrow\widetilde{\Cal V}$
\begin{Definition}\label{transverse}
A system $(\widetilde{\Cal E},\widetilde{\Cal
V},\Cal{W},\tilde{s},\Cal{N})$ is called to be a transverse Fredholm
system with index $d$ if the following hold
\bigskip

\noindent 1) For each $I\in \Cal{N}$, $\tilde{s}_I$ is a Fredholm
section, $gr(\tilde{s}_I)$ intersects transversally with the zero
section of $\tilde{E}_I\rightarrow\tilde{{\bf V}}_I$ in a
$d$-dimensional pseudomanifold $\widetilde{Z}_I$, which is also
stratified with top stratum coincides with its intersections with
the top stratum of $\tilde{{\bf V}}_I$.
\smallskip

\noindent 2) The union $Z_{\Cal{W}}=\cup_I\pi_I(\widetilde{Z}_I)$ is
compact in $\Cal{W}$. $\tilde{s}$ is called a Fredholm
multi-section.
\smallskip

\noindent 3) If the top strata of $\widetilde{Z}_I$ are oriented,
then all orientations are preserved under the partially defined
projections $\pi_J^I$.

\end{Definition}

When we generically perturb the section over the stratified spaces
we work inductively over the strata. If  $\tilde{s}_I$ is perturbed
over one stratum so that it is transverse to the zero section there,
we need to extend this perturbation to nearby strata, and hence we
need information on how that strata fit together. While McDuff
\cite{M99} showed that if there are suitable normal cones defined in
Definition \ref{cone}, for obtaining a Fredholm multi-section, it is
sufficient to perturb $gr(\tilde{s})_I$ so that the intersection of
each stratum $S'$ of $gr(\tilde{s})_I$ with a stratum $S$ of
$\tilde{{\bf V}}_I$ is transverse inside $\tilde{p}^{-1}(S)$ which
is a stratum of $\tilde{E}_I$. In application, Liu-Tian \cite{LT}
also proved that they can get a transverse Fredholm multi-section
$\bar{\partial}^{\nu}_{J,H}$ of some multi-bundle over the space of
stable $(J,H)$ maps by choosing perturbations in  finite vector
spaces $R_I$.

Then we should assemble all the local zero sets $\widetilde{Z}_I$
into a closed branched labeled pseudomanifold $Y$ which projects
onto $Z_{\Cal W}$. For each $I$, we choose a manifold with boundary
${\bf V}'_I\subset {\bf V}_I$ so that all ${\bf V}'_I$ cover $\Cal
W$. Let $\widetilde{Y}_I=\widetilde{Z}_I\cap \pi_I^{-1}({\bf
V}'_I)$. Here $\widetilde{Y}_I$ is a closed pseudomanifold and is
given the obvious first topology, and the components of
$(\widetilde{Y}_I)_{sm}$ either lie on its boundary ({\it i.e.} in
$\pi_I^{-1}(\partial {\bf V}'_I)$, or are the intersections of
strata in $\widetilde{Z}_I$ with its interior. In order for there to
be such a stratification, it suffices that the stratification of
$\Cal W$ can be refined so that each boundary $\partial ({\bf
V}'_I)$ is a union of strata. Thus the boundaries $\partial ({\bf
V}'_I)$ must intersect transversally and be in general position with
respect to the original strata of $\Cal W$. Strictly speaking, one
can arrange this for general $\Cal W$ only in the presence of
suitable normal cones. If one can find such sets $\widetilde{Y}_I$
we will say that they form a shrinking of the $\widetilde{Z}_I$.
Thus, we have stratified the $\widetilde{Y}_I$ so that their only
codimension $1$ strata lie in the boundary $\partial ({\bf V}'_I)$.
These will correspond to the branching locus of $Y$.

We then construct a topological space $Y$ such that there are
continuous maps
$$\coprod_I\tilde{Y}_I\longrightarrow Y\longrightarrow  Z_{\Cal W},$$
where $Y$ consists of points $[y]$ which are the equivalence classes
under the equivalence relation given by
$$\tilde{x}_I\sim\tilde{y}_J,\ \ \ \ {\rm if}\ J\subset I,\ \ \pi_J^I(\tilde{x}_I)=\tilde{y}_J.$$
The first topology on $Y$ is the quotient topology, and the strata
of its second (smooth) topology are formed by the images of strata
in the $\tilde{Y}_I$, subdivided if necessary. It will be convenient
to refine the stratifications of the $\tilde{Y}_I$ so that the
projections $q_I : \tilde{Y}_I\rightarrow Y$ take strata to strata.
This introduces new codimension 1 strata in the $\tilde{Y}_I$ coming
from the boundaries of the $\tilde{Y}_K$ for $I\subset K$.

In order to define a cycle, we still need a suitable labeling of the
top components, {\it i.e.} to define a positive rational labeling
function. For a top stratum $S$ lying the image of $q_I :
\tilde{Y}_I\rightarrow Y$, we define
$$\lambda_I([y])=\frac{|\{q_I^{-1}([y])\}|}{|\Gamma_I|},\ \ [y]\in S,$$
where $|\cdot|$ denotes the cardinality of a set or group. McDuff
proved the following result, for the reader's convenience we restate
her proof.
\begin{Proposition}
Under the setting as above the labeling function $\lambda_I$
descends to a function $\lambda$ on $Y$, and $(Y,\lambda)$ can be
regarded as a compact labeled branched pseudomanifold as is defined
in Definition \ref{label}.
\end{Proposition}
Proof. Since for $J\subset I$ and $y\in q_I(\tilde{Y}_I)\cap
q_J(\tilde{Y}_J)$, $\tilde{s}_I$ is the pull-back of $\tilde{s}_J$
over $(\pi_J^I)^{-1}(\tilde{\bf V}_J)$ such that
$$|\{q_I^{-1}([y])\}|=|\{q_J^{-1}([y])\}|\cdot|\Gamma_{I-J}|,$$
and $|\Gamma_I|=|\Gamma_J|\cdot|\Gamma_{I-J}|$, we see that
$\lambda_I([y])=\lambda_J([y])$, thus we get a function $\lambda$ on
$Y$.

To check the branching condition in Definition \ref{label} on the
codimension 1 components of $Y$, we have to understand how the
generating equivalence $\tilde{x}_I\sim\tilde{y}_J$, $J\subset I$,
effect on such condition. Note first that a class $[y]$ lies in a
codimension 1 stratum $B_l$ in $Y$, if for one $K$ the
representatives of $[y]$ in $\tilde{Y}_K$ lie on the boundary
stratum $B_K$ and if for all other $J$ the representatives in
$\tilde{Y}_J$ lie in the interior of $\tilde{Y}_J$. We call the side
of $B_K$ that meets the interior of $\tilde{Y}_K$ the positive side,
and denote by $\{S^+_{K,\alpha}\}$ for the set of all top strata in
$\tilde{Y}_K$ whose closures contain representative of $[y]$. For
all other $J$, we denote by $\{S^{\pm}_{J,\alpha}\}$ for the set of
all top strata in $\tilde{Y}_J$ containing representative of $[y]$,
with sign assigned in the obvious way. By the Lemma \ref{VI}, the
intersection $I_y=J\cap K$ is nonempty and everything is essentially
pulled back from $\tilde{V}_{I_y}$.  For instance, two strata in
$\{S^{+}_{J,\alpha}\}$ are identified in $Y$ if and only if they
have the same image in  $\tilde{V}_{I_y}$. Therefore, we only
consider the effect of inclusions of the form $I_y\subset I$ on the
branching condition. We distinguish three cases: i) $K\neq I_y$ or
$I$; ii) $K=I$; iii) $K=I_y$.

For case i), assume $J$ is either $I_y$ or $I$, there is a bijective
correspondence between the components $\{S^{+}_{J,\alpha}\}$ and
$\{S^{-}_{J,\alpha}\}$ on both sides that commutes with the
identifications coming from the equivalence $\sim$. So the labelings
are the same on both sides.

For case ii), each stratum in $\{S^{+}_{K,\alpha}\}$ is mapped
bijectively by $\pi_{i_Y}^{k}$ onto a stratum in
$\{S^{+}_{I_y,\alpha}\}$, and each of the latter strata is covered
exactly $|\Gamma_{K-I_y}|$ times. Also both sides of $B$ in $Y$ are
the same and no real branching occurs either.

For case iii), the situation is different. The components in
$\{S^{+}_{I,\alpha}\}$ are identified via $\pi_K^I$ with components
in $\{S^{+}_{K,\alpha}\}$, while no identification is put on the
components $\{S^{+}_{I,\alpha}\}$ on the other side. Since $\pi_K^I$
is $|\Gamma_{I-K}|$ to 1, $|\Gamma_{I-K}|$ components on the
negative side of $B$ in $Y$ will correspond to each component on the
positive side. One can check that the sum of labels on both sides is
the same.  \qed

Recall that such a $d$-dimensional branched pseudomanifold $Y$ can
be considered as a cycle in the sense that every partially smooth
map $f$ from $Y$ to a closed manifold $X$ defines a rational class
$f_*([Y])\in H_d(X)$, say {\it virtual cycle} (c.f. Lemma \ref{5.1}
or \cite{M99}).

\section{Transversality.}
Let us come back to our problem. In order to apply arguments showed
above to our case of moduli space of stable connecting orbits, we
need firstly construct a transverse Fredhlom system from
$P\Cal{M}(\tilde{x}_-,\tilde{x}_+)\subset
\Cal{B}(\tilde{x}_-,\tilde{x}_+)$ as defined in the Definition
\ref{transverse}. Because of the possibility of bubbling-off
multiple covered holomorphic spheres with negative first Chern class
in a connecting orbit $u\in \Cal{M}(\tilde{x}_-,\tilde{x}_+)$, even
for generic pair $(J,H)$, the linearized map $D_u$ in general does
not be surjective. However, one can show that under a locally
non-equivariant perturbation of the $\bar{\partial}_{J,H}$-operator
constricted in a small neighborhood $U_\epsilon$ in the ambient
space $\Cal{B}(\tilde{x}_-,\tilde{x}_+)$, denoted by
$\bar{\partial}_{J,H}+\nu$, the regularity will hold. For instance,
for defining the classical Floer homology for some chain complex
generated by nondegenerate periodic solutions of Hamiltonian flow,
Liu-Tian (c.f. section 3 in \cite{LT}) showed the full details to
realize the regularization by non-equivariant perturbations, and
define the local virtual moduli space ${\Cal M}_\epsilon^\nu$ which
has the expected dimension. Since the only difference between in our
case and in the classical one is that the boundary condition
(\ref{002}) is related to a symplectomorphism $\phi$, which may not
be identity, {\it i.e.} the solutions $x_\pm$ are unnecessarily
periodic, and this difference has little effect on all the process
of construction, we then only give a sketchy and suggestive
description below.

For any two $\tilde{x}_-,\tilde{x}_+\in Crit(F)$, we suppose that
the unparameterized stable orbit $[V]$ connecting $\tilde{x}_-$ and
$\tilde{x}_+$.  Recall $(\widetilde{U}_{\epsilon}^D(V,{\bf
H}),\Gamma_V,\pi_V)$ is a local uniformizer for $[V]$ in ${\Cal
W}^D$, and  denote $\widetilde{U}_{\epsilon}(V,{\bf H})=\cup_D
\widetilde{U}_{\epsilon}^D(V,{\bf H})$ with automorphism group
$\Gamma_V$, and we saw that $U_{\epsilon}(V,{\bf
H})=\widetilde{U}_{\epsilon}(V,{\bf H})/\Gamma_V$ can be regarded as
a small neighborhood of $[V]$ in ${\Cal W}=\bigcup_{[V]\in P{\Cal
M}}U_{\epsilon_V}(V,{\bf H})$ which is a covering of
$P\Cal{M}(\tilde{x}_-,\tilde{x}_+)$.

We then define locally the bundles $\tilde{E}^D(V)$ and
$\tilde{E}(V)$ over $\widetilde{U}_{\epsilon}^D(V,{\bf H})$ and
$\widetilde{U}_{\epsilon}(V,{\bf H})$ as follows. For each
$\hat{V}\in\widetilde{U}_{\epsilon}^D(V,{\bf H})$ or
$\widetilde{U}_{\epsilon}(V,{\bf H})$, the fiber is
$$\Cal{F}_{\hat{V}}=\{\xi|\xi\in
L^{k-1,p}(\Lambda^{0,1}(\hat{V}^*TM))\},$$ where the $L^{k-1,p}$
norm is suitably introduced by gluing process from the metric on the
domain. We can see that for fixed $D$, $\tilde{E}^D(V)$ is a locally
trivial Banach bundle over $\widetilde{U}_{\epsilon}^D(V,{\bf H})$.
So the topology of $\tilde{E}(V)$ is well-defined when restricted to
each stratum of $\widetilde{U}_{\epsilon}^D(V,{\bf H})$. As we
claimed before, the $\bar{\partial}_{J,H}$-operator can be
considered as a $\Gamma_V$-equivariant section of the bundle
$\tilde{E}(V)\rightarrow \widetilde{U}_{\epsilon}(V,{\bf H})$, it is
smooth on each stratum $\widetilde{U}_{\epsilon}^D(V,{\bf H})$. The
zero sets $\bar{\partial}^{-1}_{J,H}(0)$ in
$\widetilde{U}_{\epsilon}^D(V,{\bf H})$ and
$\widetilde{U}_{\epsilon}(V,{\bf H})$ are projected to
$$P\Cal{M}^D(\tilde{x}_-,\tilde{x}_+)\cap U_{\epsilon}^D(V,{\bf
H})\subset {\Cal B}^D(\tilde{x}_-,\tilde{x}_+)$$ and
$$P\Cal{M}(\tilde{x}_-,\tilde{x}_+)\cap U_{\epsilon}(V,{\bf
H})\subset {\Cal B}(\tilde{x}_-,\tilde{x}_+),$$ respectively.

We first consider those $[V]$ with fixed intersection pattern $D$.
We can define a coordinate chart of
$\widetilde{U}_{\epsilon}^D(V,{\bf H})$ and a trivialization of
$\tilde{E}^D(V)$. First, for fixed parameter $\alpha$ with
intersection pattern $D$, we assume $x_{l,j}$ is the added marked
point of a component $\Sigma^l_{\alpha}$ of $\Sigma_{\alpha}$, and
we denote the tangent space of ${\bf H}_{l,j}$ at $x_{l,j}$ by
$H_{l,j}$. Recall the fact that we are restricting to a slice for
the action of reparametrization group, we define the space
$$L^{k,p}(V_{\alpha}^*TM,H)=\{\xi=\xi_{\alpha}|\xi\in L^{k,p}(V_{\alpha}^*TM)
,\ \xi_l(x_{l,j})\in H_{l,j}\},$$ and the set
$$W_{\epsilon}^{\alpha}=\{\xi\in L^{k,p}(V_{\alpha}^*TM,H),\ \|\xi\|_{k,p}<\epsilon\}.$$
It is a smooth coordinate chart for
$\widetilde{U}_{\epsilon_1}^{\alpha}(V_{\alpha},{\bf H})$ near
$V_{\alpha}$ for $\epsilon\ll\epsilon_1$. So
$W_\epsilon=\cup_{\alpha\in D}W_{\epsilon}^{\alpha}$ is the
coordinate chart for $\widetilde{U}_{\epsilon_1}^D(V,{\bf H})$. In
fact, if we denote
$$\Lambda_\epsilon=\{\alpha|\ \alpha\in D, \|\alpha\|<\epsilon\},$$
the $W_\epsilon$ has a splitting
$W_\epsilon^{\alpha=0}\times\Lambda_\epsilon$ which can be regarded
as the local coordinate of $\widetilde{U}_{\epsilon_1}^D(V,{\bf
H})$.

As the standard method (c.f. \cite{MS}\cite{LT}), we can get a
trivialization of the bundle $\tilde{E}^{\alpha}(V)$ and
$\tilde{E}^D(V)$, we denote the trivialization by
$$\gamma_D:\ \widetilde{U}_{\epsilon}^D(V,{\bf H})\times L^{k-1,p}(\Lambda^{0,1}(V^*TM))\rightarrow\tilde{E}^D(V).$$

Under these local charts and trivialization, for fixed $\alpha$, we
can write the section $$\bar{\partial}_{J,H}:\
\widetilde{U}_{\epsilon}^{\alpha}(V,{\bf
H})\rightarrow\tilde{E}^{\alpha}(V)$$ as a nonlinear map
$$S_{\alpha}:\ W_\epsilon^{0}\rightarrow L^{k-1,p}(\Lambda^{0,1}(V^*TM)).$$
Then the problem of transversality is to consider whether its
linearized map is surjective. In particular, if $\epsilon$ is
sufficiently small, we only need to componentwise deal with the
operator
$$T_l=DS_{\alpha=0}(0):\ L^{k,p}(V_l^*TM,H)\rightarrow
L^{k-1,p}(\Lambda^{0,1}(V_l^*TM)),$$ which is a linear elliptic
operator, it is also a Fredholm operator which can be proved via
standard arguments (c.f.\cite{MS}). While it is not surjective in
general duo to the appearing of multiple covered bubble sphere with
negative first Chern class. However, the cokernel $K_l=K_l(V)$ of
$T_l$ is finite dimensional.\footnote{For later extending $K_l$ to a
vector bundle over $\widetilde{U}_{\epsilon}(V,{\bf H})$, we may
assume its vectors vanish at each bubble point (c.f.\cite{LT}).} Let
$K=\oplus_{l=1}^LK_l$, then $L^{k-1,p}(\Lambda^{0,1}(V^*TM))=K\oplus
Im(T)$. On each main component, for generic pair $(J,H)$, $T_m$ is
surjective, so we only deal with bubble components. Then Liu-Tian in
\cite{LT} showed the method of enlarging the domain of $T$ to
realize surjectivity. Roughly speaking, we can define a finite
dimensional vector space $R\subset L^{k-1,p}(\Lambda^{0,1}(V^*TM))$
such that the new linear map
$$T\oplus{\bf I}:\
L^{k,p}(V_l^*TM,H)\oplus R\rightarrow
L^{k-1,p}(\Lambda^{0,1}(V_l^*TM))$$ is surjective and the kernel of
it is the same as the kernel of $T$, where ${\bf I}:\ R\rightarrow
L^{k-1,p}(\Lambda^{0,1}(V_l^*TM))$ is inclusion.

Then we can extend $R$ over $\widetilde{U}_{\epsilon}^{0}(V,{\bf
H})$ and $\widetilde{U}_{\epsilon}^D(V,{\bf H})$ with $D=D(f)$. Thus
when $\|\alpha\|$ are small enough, we can get $R(V)\subset
L^{k-1,p}(\Lambda^{0,1}(V_{\alpha}^*TM))$. If dim$R=r$, then we
obtained a $r$-dimensional vector bundle $R$ over
$\widetilde{U}_{\epsilon}^D(V,{\bf H})$. Also we have a surjective
linear map $$T_{\alpha}\oplus{\bf I}_{\alpha}:\
L^{k,p}(V_{\alpha}^*TM,H)\oplus R(V_{\alpha})\rightarrow
L^{k-1,p}(\Lambda^{0,1}(V_{\alpha}^*TM)).$$ Consequently, by
implicit function theorem we know that the following moduli space
$$\widetilde{\Cal M}^D_{R,\epsilon}(\tilde{x}_-,\tilde{x}_+)=\{\hat{V}|\ \hat{V}\in
\widetilde{U}^D_{\epsilon}(V,{\bf H}),\
\bar{\partial}_{J,H}\hat{V}\in R\}$$ is a smooth manifold of
dimension $r+\mu_{rel}(\tilde{x}_-,\tilde{x}_+)-1$.

\bigskip
When we consider the changes of the topological type of the domains,
we use the gluing parameter $(t,\tau)$ introduced in subsection 4.2.
Similarly, for a fixed parameter $(\alpha,t,\tau)$, we can define
$W_\epsilon^{(\alpha, t,\tau)}$ as a coordinate chart of
$\widetilde{U}_{\epsilon}^{(\alpha, t,\tau)}(V,{\bf H})$, the
trivialization of the bundle $\tilde{E}^{(\alpha,
t,\tau)}\rightarrow\widetilde{U}_{\epsilon}^{(\alpha,
t,\tau)}(V,{\bf H})$, and under the coordinate and trivialization we
can regard the $\bar{\partial}_{J,H}$ section of this bundle as a
nonlinear map
$$S_{(\alpha, t,\tau)}:\
W_\epsilon^{(\alpha, t,\tau)}\rightarrow
L^{k-1,p}(\Lambda^{0,1}(V_{(\alpha, t,\tau)}^*TM)),$$ whose
linearized map
$$T_{(\alpha, t,\tau)}=DS_{(\alpha, t,\tau)}(0):\
L^{k,p}(V_{(\alpha, t,\tau)}^*TM,H)\rightarrow
L^{k-1,p}(\Lambda^{0,1}(V_{(\alpha, t,\tau)}^*TM))$$ is a Fredholm
operator. And we can define a finite dimensional subspace
$R_{(\alpha, t,\tau)}=R(V_{(\alpha, t,\tau)})\subset
L^{k-1,p}(\Lambda^{0,1}(V_{(\alpha, t,\tau)}^*TM))$ and embedding
${\bf I}$ such that the operator $$T\oplus {\bf I}:\
L^{k,p}(V_{(\alpha, t,\tau)}^*TM,H)\oplus R_{(\alpha,
t,\tau)}\rightarrow L^{k-1,p}(\Lambda^{0,1}(V_{(\alpha,
t,\tau)}^*TM))$$

Then the critical thing is to show that when the parameters change
in a small neighborhood
$$\Lambda_{\delta}=\{(\alpha, t,\tau)|\ \|(\alpha, t,\tau)\|<\delta \},$$ one can still find
$R_{(\alpha, t,\tau)}$ and a partially smooth family of embedding
${\bf I}_{(\alpha,t,\tau)}$ such that the linearized operator
$$T_{(\alpha,t,\tau)}\oplus {\bf
I}_{(\alpha,t,\tau)}:\ L^{k,p}(V_{(\alpha,t,\tau)}^*TM,H)\oplus\
R_{(\alpha,t,\tau)}\rightarrow
L^{k-1,p}(\Lambda^{0,1}(V_{(\alpha,t,\tau)}^*TM))$$ is surjective,
{\it i.e.} with uniform estimates for its right inverse as
$(\alpha,t,\tau)$ varies in the sufficiently neighborhood
$\Lambda_{\delta}$.

In order to get the desired uniform estimate, Liu-Tian (c.f. section
3 in \cite{LT}) used some exponential weighted equivalent norms
$\|\cdot\|_{\chi;k,p}$ on $L^{k,p}(V_{(\alpha,t,\tau)}^*TM,H)$ and
$L^{k-1,p}(\Lambda^{0,1}(V_{(\alpha,t,\tau)}^*TM))$. The same
argument can apply to our case, and we can similarly prove that
there exists a right inverse $G$ of $T_{(\alpha,t,\tau)}\oplus {\bf
I}_{(\alpha,t,\tau)}$ and a constant $c=c(V)$ depending only on $V$
such that for sufficiently small $\delta$,
$(\alpha,t,\tau)\in\Lambda_{\delta}$, we have
$$\|G_{(\alpha,t,\tau)}\xi\|_{\chi;k,p}<c(V)\|\xi\|_{\chi;k,p}$$
for $\forall\ \xi\in L^{k-1,p}(\Lambda^{0,1}(V_{(\alpha,t,\tau)}))$.
That is equivalent to say, if we shrink
$\widetilde{U}_\epsilon(V,{\bf H})$ to be a sufficiently small
neighborhood, then there exists a uniformly bounded family of right
inverses to $T_{\hat{V}}\oplus {\bf I}_{\hat{V}}$ as $\hat{V}$
varies in $\widetilde{U}_\epsilon(V,{\bf H})$. And in this small
neighborhood $\widetilde{U}_\epsilon(V,{\bf H})$, we can identify
$R_{(\alpha,t,\tau)}$ with $R_V=R_{(0,0,0)}$.

Then we can define $${\Cal O}={\Cal
O}_{\widetilde{U}_\epsilon(V,{\bf H})}=\widetilde{U}_\epsilon(V,{\bf
H})\times R_V$$ and a projection $$P:\ {\Cal
O}\rightarrow\widetilde{U}_\epsilon(V,{\bf H}),$$ for the pullback
bundle $P^*(\tilde{E}_V)\rightarrow\widetilde{U}_\epsilon(V,{\bf
H})\times R_V$, we can construct a section as\footnote{Here for
simplicity we admit to abuse the notations. Actually, the first
summand is the $\bar{\partial}_{J,H}$-operator for each main
component and is the $\bar{\partial}_J$-operator for each bubble
component.}
\begin{equation}\label{section}
s(\hat{V},\nu)=\bar{\partial}_{J,H}(\hat{V})+{\bf I}_{\hat{V}}(\nu),
\end{equation}
where $\nu\in R_V$. From the construction above, its linearized
operator is surjective at all points $(\hat{V},\nu)\in {\Cal O}$.
Then by using some variant of gluing as in subsection 3.3 of
\cite{LT} one can see that the zero set $s^{-1}(0)$ is a open
partially smooth pseudomanifold of dimension
$r+\mu_{rel}(\tilde{x}_-,\tilde{x}_+)$, whose components are the
intersections of $s^{-1}(0)$ with each stratum
$\widetilde{U}^D_\epsilon(V,{\bf H})$. So we obtain a new transverse
Fredholm system $(P^*(\tilde{E}_V),{\Cal O},s)$. Since $R_V$ is
finite dimensional, the Sard-Smale theorem (c.f. \cite{MS}) says
that for generic $\nu\in R_V$ we have a transverse Fredholm section
$$s^{\nu}:\ \tilde{E}_V\rightarrow\widetilde{U}_\epsilon(V,{\bf H}),$$
satisfying $s=P^*(s^{\nu})$. Then we denote its zero set by
$\widetilde{Z}_V^{\nu}=(s^{\nu})^{-1}(0)$ and so we have

\begin{Proposition}\label{Z-nu}
For generic $\nu\in R_V$, the section $s^{\nu}$ is transverse
Fredholm and its zero set $\widetilde{Z}_V^{\nu}$ has the structure
of an open pseudomanifold of dimension
$d=\mu_{rel}(\tilde{x}_-,\tilde{x}_+)-1$.
\end{Proposition}

Now let
$$\widetilde{\Cal M}^{\bar{D}_1}_{R,\epsilon}(\tilde{x}_-,\tilde{x}_+)=
\bigcup_{(\alpha,t,\tau)\in D',D\le D'\le D_1}\widetilde{\Cal
M}^{(\alpha,t,\tau)}_{R,\epsilon}(\tilde{x}_-,\tilde{x}_+).$$ If we
denote $n_t$ and $n_\tau$ for the numbers of zero components of the
gluing parameter $t$ and $\tau$ for a generic
$\bar{D}_1\in\Lambda_{\bar{D}_1}$, where
$\bar{D}_1=\{(\alpha,t,\tau)|\ (\alpha,t,\tau)\in D',\ D\le D'\le
D_1\}$, then with similar argument as above or, with more details as
Liu-Tian \cite{LT} proved, we see that for generic choice $\nu\in
R$, the moduli space of unparameterized stable
$(J,H,\nu)$-connecting orbits
$$\widetilde{\Cal
M}^{\bar{D}_1,\nu}_{R,\epsilon}(\tilde{x}_-,\tilde{x}_+)\subset\widetilde{\Cal
M}^{\bar{D}_1}_{R,\epsilon}(\tilde{x}_-,\tilde{x}_+)$$ is a cornered
partially smooth pseudomanifold with correct dimension
$$\mu_{rel}(\tilde{x}_-,\tilde{x}_+)-1-\Sigma(2n_t+n_\tau).$$
Moreover, the transversality can be achieved for all $D'$ with $D\le
D'\le D_1$ simultaneously.

\section{Virtual moduli cycle and Floer homology. }

\subsection{Constructing virtual moduli cycle}

Since $P\Cal{M}(\tilde{x}_-,\tilde{x}_+)$ is compact, we can take a
finite union of the covering of $\Cal W$ as
$\{U_i=U_{\epsilon_{V_i}}(V_i,{\bf H}), i=1,\cdots, w\}$, and we use
$\tilde{U}_i$ to denote its uniformizer with covering group
$\Gamma_i$.

Now just as the construction in subsection 4.3, we can define
$U_I=\cap_{i\in I}U_i$ where $I=\{1,\cdots,w\}$, and the fiber
product $\tilde{U}_I$ as in Definition \ref{product}, so we can get
suitable ${\bf V}_I\subset U_I$ and $\tilde{{\bf V}}_I$ with
quotient map $\pi:\tilde{{\bf V}}_I\rightarrow {\bf
V}_I\simeq\tilde{{\bf V}}_I/\Gamma_I$ and the projection
$\pi_J^I:\tilde{{\bf V}}_I\rightarrow \tilde{{\bf V}}_J$ for
$J\subset I$. Thus as in Definition \ref{m-fold} we have a
multi-fold atlas $\widetilde{\Cal V}$ for $\Cal{W}$.

Then as above we have a local orbifold bundle $E_i$ over $U_i$. For
each $[V]\in U_i$, the fiber ${\Cal F}_i|_{[V]}$ over $[V]$ consists
of all elements of $L^{k-1,p}(\Lambda^{0,1}(V^*TM)),\ V\in[V]$
modulo equivalence relation induced by pull-back of sections coming
from identificaion of the domains of $V_i$, where
$\Lambda^{0,1}(V^*TM))$ is the bundle of $(0,1)$-forms on $\Sigma$
with respect to the complex structure on $\Sigma$ and the given
compatible almost complex structure $J$ on $(M,\omega)$. Then the
local uniformizer $\tilde{E}_i$ of $E_i$ is given by the union of
$L^{k-1,p}(\Lambda^{0,1}(\tilde{V}_i^*TM))$, $\tilde{V}_i\in
\tilde{U}_i$. The $\Gamma_i$ also acts on $\tilde{E}_i$ so that
$E_i=\tilde{E}_i/\Gamma_i$. In this way, we can reinterpret the
$\bar{\partial}_{J,H}$-operator as a collection of
$\Gamma_i$-equivalent sections of these local orbifold bundles
$(E_i,U_i)$.

Then as in Lemma \ref{VI} we can choose subcovers $\{U_i^0\}$ and
$\{{\bf V}_I\}$ and a suitable partition of unity $\beta_i$ on $\Cal
W$ corresponding to the covering $\{U_i^0\}$. Let ${\Cal
E}=\cup_iE_i$, by the construction of section 4, we obtain a
multi-bundle $\tilde{p}:\ \widetilde{\Cal
E}\rightarrow\widetilde{\Cal V}$.

Denote $R_i=R_{V_i}$, and the projection to the first factor
$P:\tilde{U}_i\times R_i\rightarrow\tilde{U}_i$. We can define the
section of the pullback bundle $P^*(\tilde{E}_i)\rightarrow
\tilde{U}_i\times R_i$ as
$$\iota(i)(\tilde{V},\nu_i)=\beta_i([V])\cdot{\bf
I}_{V}(\nu_i).$$ Since there are only finite small neighborhoods, we
can choose a sufficiently small $\epsilon$ such that for all $i$ and
generic $\nu_i$, $\bar{\partial}_{J,H}+\iota(i)$ are transverse
Fredholm sections of the corresponding pullback bundles. Now we set
$R=\oplus_{i=1}^w R_i$ and choose its a small subset
$R_\epsilon=\{\nu\in R| |\nu_i|\le\epsilon, \ \forall i\}$. Then let
${\Cal W}_\epsilon={\Cal W}\times R_\epsilon$, we have a
corresponding multi-fold atlas $\widetilde{\Cal V}_\epsilon$. Then
we can give a multi-bundle structure to $P^*(\widetilde{\Cal
E})\rightarrow\widetilde{\Cal V}_\epsilon$. Since $\widetilde{\Cal
V}_\epsilon$ is a good cover, the compatibility condition
(\ref{compatibility}) holds, we get a multi-section $\tilde{s}$ of
this multi-bundle as
\begin{equation}
\tilde{s}(\tilde{V},\nu)=\bar{\partial}_{J,H}(\tilde{V})+\sum_i\iota(i)(\tilde{V},\nu_i).
\end{equation}
It is easy to see that $(\widetilde{\Cal E},\widetilde{\Cal
V}_\epsilon,\Cal{W}_\epsilon,\tilde{s})$ is a transverse Fredholm
system with index $d=r+\mu_{rel}(\tilde{x}_-,\tilde{x}_+)$ as
defined in Definition \ref{transverse}. And it follows that when
$\epsilon$ is small enough for a generic choice of the perturbation
$\nu\in R_{\epsilon}$ the section $\bar{\partial}_{J,H}+\tilde{\nu}$
of $\widetilde{\Cal E}\rightarrow\widetilde{\Cal V}$ is Fredhom and
have the same zero set as $\tilde{s}$. Consequently, corresponding
to Proposition \ref{Z-nu}, and by using the method of constructing
the branched labeled pseudomanifold with boundary $Y$ from the zero
set of the multi-section $\tilde{s}$ in subsection \ref{construct}
we have
\begin{Proposition}\label{MMM-nu}
For generic $\nu\in R_V$, the multi-section $\tilde{s}^{\nu}$
satisfying $\tilde{s}=P^*(\tilde{s}^{\nu})$ is transverse Fredholm
and its zero sets $\widetilde{Z}_I^{\nu}$ fit together to give a
compact branched and labeled pseudomanifold $Y=P{\Cal
M}^\nu(\tilde{x}_-,\tilde{x}_+)$ with boundary
$B^\nu(\tilde{x}_-,\tilde{x}_+)$, which is a relative virtual moduli
cycle of dimension $d=\mu_{rel}(\tilde{x}_-,\tilde{x}_+)-1$.
\end{Proposition}
Proof. The conclusion is straightforward. We only show an explicit
proof by Liu-Tian for the compactness of $Y=P{\Cal
M}^\nu(\tilde{x}_-,\tilde{x}_+)$.  From the construction in
subsection \ref{construct} we know that $P{\Cal M}^\nu$ is projected
onto $Z_{\Cal W}=\overline{P{\Cal M}}^\nu$ in ${\Cal W}$. Since
$$P{\Cal M}(\tilde{x}_-,\tilde{x}_+)\subset\cup_{i=1}^wU_i^0,$$
where the $U_i^0$ is constructed in the Lemma \ref{VI}, we see that
$\bar{\partial}_{J,H}$ never becomes zero along the boundary
$$\partial(\cup_{i=1}^wU_i^0)=\overline{\cup_{i=1}^wU_i^0}\setminus\cup_{i=1}^wU_i^0.$$
But $\nu\equiv 0$ along $\partial(\cup_{i=1}^wU_i^0)$, so we have
$$\overline{P{\Cal M}}^\nu\subset\cup_{i=1}^wU_i^0).$$ Let
$\{\bar{V}_i\}_{i=1}^\infty$ be a sequence of $\overline{P{\Cal
M}}^\nu$. We may assume that all $\bar{V}_i$ are contained in ${\bf
V}_I$ for some $I\in{\Cal N}$. If we can show that for the
corresponding sequence $\{V_i\}_{i=1}^\infty$, with
$\bar{V}_i=\pi^I(V_i)$, in $P{\Cal M}^\nu_I=P{\Cal
M}^\nu\cap\tilde{\bf V}_I$, all sections $\nu_I$ of the bundle
$\tilde{E}_I\rightarrow\tilde{\bf V}_I$ has a uniform bounded
$W^{k,p}$-norm, then by the basic elliptic techniques, there exists
a $V_\infty\in\tilde{U}_I$ such that some subsequence of
$\{V_i\}_{i=1}^\infty$ is weakly $C^\infty$-convergent to
$V_\infty$. Then $\overline{P{\Cal M}}^\nu$ is compact.

Since $\nu=\sum_{i,j}a_{ij}e_{ij}$ with $\{e_{ij};j=1,\cdots,n_i\}$
being the basis of $R_i$, and $\nu\in R_\epsilon$, we see that
$|a_{ij}|$ are bounded. So we only need to prove that all the
lifting $\{e_{ij}\}_I$ over $\tilde{\bf V}_I$ of $e_{ij}$, which is
defined over $\tilde{U}_i$ originally, are still bounded. We still
only consider the case $i\in_I$. Actually, if we can prove that all
changes of coordinates between $\tilde{U}_i$'s are induced from
those reparametrizations that stay inside a compact subset of
$G_\Sigma$, then the boundedness of $\|\nu_I\|_{k,p}$ will follow.

Now we denote the closure of $U_i$ by $U^c_i$, $i=1,\cdots,w$, and
$U_{ij}^c=U^c_i\cap U^c_j$. Let $\tilde{U}^c_i$ and
$\tilde{U}_{\bar{i}j}^c\subset\tilde{U}^c_i$ be the lifting of them
in the uniformizer $\tilde{U}_i$. Then we denote the compact set
$$P{\Cal M}_{ij}=P{\Cal M}(\tilde{x}_-,\tilde{x}_+)\cap U_{ij}^c.$$
Let
$$\{Z_k^{ij}|\ Z_k^{ij}\subset{\Cal B}(\tilde{x}_-,\tilde{x}_+),\ k=1,\cdots,m^{ij}  \}$$
be an open covering of $P{\Cal M}_{ij}$ in ${\Cal
B}(\tilde{x}_-,\tilde{x}_+)$ such that each component of
$\pi_i^{-1}(Z_k^{ij})$ and $\pi_j^{-1}(Z_k^{ij})$ in $\tilde{U}_i^c$
and $\tilde{U}_j^c$ is a uniformizer of $Z_k^{ij}$, respectively.
Now for each fixed pair of components of $\pi_i^{-1}(Z_k^{ij})$ and
$\pi_j^{-1}(Z_k^{ij})$, the equivalence between them are induced by
some automorphisms of domain which are contained in a compact subset
of $\prod G_\Sigma$. Then let $Z^{ij}=\cup_kZ_k^{ij}$ and for each
$i=1,\cdots,w$, replace $U_i$ by $$(U_i\setminus\cup_{k\neq i}
U^c_{ik})\bigcup_{k\neq i}Z^{ik}.$$ They still form an open covering
of $P{\Cal M}(\tilde{x}_-,\tilde{x}_+)$ and all previous
constructions still work. Now all changes of coordinates are induced
by a compact subset. \qed

Recall that each top stratum $S$ of $Y$ lying in the image
$q_I:\widetilde{Y}_I\rightarrow Y$ is associated a positive rational
label $\lambda_I$. They can fit together to give each
$\delta_i$-oriented top component $M_i$ a label $\lambda_i$. Then we
can define a rational number for the virtual moduli space
$$\#(Y)=\#(P{\Cal M}^{\nu})=\sum_i\delta_i\lambda_i.$$

We say any compact branched and labeled pseudomanifold $P{\Cal
M}^\nu(\tilde{x}_-,\tilde{x}_+)$ as above constructed is the {\it
regularization} or {\it virtual moduli space} of the stable moduli
space $P{\Cal M}(\tilde{x}_-,\tilde{x}_+)$. In particular, we care
about such $d=\mu_{rel}(\tilde{x}_-,\tilde{x}_+)-1=0 $ and 1
dimensional pseudomanifolds
$$P{\Cal
M}^\nu(\tilde{x}_i,\tilde{x}_{i+1})\ \ \ {\rm and}\ \ \ P{\Cal
M}^\nu(\tilde{x}_i,\tilde{x}_{i+2}).$$ The former is obviously a
finite set. And the latter is a branched and labeled
1-pseudomanifold with boundary which, via Floer's gluing method,
consists of pairs $[V\#U]$ with $[V]\in P{\Cal
M}^\nu(\tilde{x}_i,\tilde{x}_{i+1})$ and $[U]\in P{\Cal
M}^\nu(\tilde{x}_{i+1},\tilde{x}_{i+2})$. Recall that in section 5
we can associate to each boundary point a rational number
$\rho([V\#U])$, then from the Lemma \ref{5.2} we see that the total
oriented number of its boundary is
$\#(B^\nu(\tilde{x}_{i},\tilde{x}_{i+2}))=\sum_{x\in
B^{\nu}}\rho(x)=0$. Thus, we naturally have the following result
\begin{Corollary}\label{COR}
$1^\circ$. If the relative index
$$\mu_{rel}(\tilde{x}_-,\tilde{x}_+)=\mu_{CZ}(\tilde{x}_+)-\mu_{CZ}(\tilde{x}_-)=1,$$
then $P{\Cal M}^\nu(\tilde{x}_-,\tilde{x}_+)$ is a finite set;\\
$2^\circ$. If $\mu_{rel}(\tilde{x}_-,\tilde{x}_+)=2$, then in the
sense of partially smooth category, the oriented boundary
$B^\nu(\tilde{x}_-,\tilde{x}_+)=\partial(P{\Cal
M}^\nu(\tilde{x}_-,\tilde{x}_+))$ is a finite set with the total
oriented number $\#(B^\nu(\tilde{x}_-,\tilde{x}_+))=0$. Moreover,
$$\partial(P{\Cal
M}^\nu(\tilde{x}_-,\tilde{x}_+))=\sum_{\mu_{rel}(\tilde{x}_-,\tilde{y})=1;\mu_{rel}(\tilde{y},\tilde{x}_+)=1}
P{\Cal M}^\nu(\tilde{x}_-,\tilde{y})\times P{\Cal
M}^\nu(\tilde{y},\tilde{x}_+).$$

\end{Corollary}
Remark. This Corollary can be generalized to the case for
$t$-dependent pair $(J_t,H_t)$ that will be used in the continuation
argument.

\subsection{Define the Floer-type homology}

We come to define the Floer chain complex. As the usual way, we
first define a graded $\Q$-space
$C_*=C_*(J,H,\phi)=\oplus_nC_n(J,H,\phi)$ as follows. Recall that
$F=F_H$ is the functional defined in section 2 and $Crit(F)$ is the
set of all critical points. We denote by $Crit_n(F)$ the subset of
all $\tilde{x}\in Crit(F)$ with Conley-Zehnder index
$\mu_{CZ}(\tilde{x})=n$. Let $C_n=C_n(J,H,\phi)$ be the set of all
formal sums $$\xi=\sum_{\tilde{x}\in
Crit_n(F)}\xi_{\tilde{x}}\tilde{x},$$ where $\xi_{\tilde{x}}\in\Q$
such that for any constant $c>0$ $$\#\{\tilde{x}\in Crit_n(F)|\
\xi_{\tilde{x}}\neq 0,\ F(\tilde{x})\geq c \}<\infty.$$  Thus $C_*$
is an infinite dimensional vector space over $\Q$ in general,
however, if we introduce the so-called Novikov ring
$\Lambda_{\omega,\phi}$ (which is a field here) as follows, we will
see that it is a finite dimensional vector space over
$\Lambda_{\omega,\phi}$.

Recall that $\Gamma$ is the covering group introduced in section 2,
and we assume that there is an injective homomorphism
$i:\Gamma\rightarrow\pi_2(M)$. Then the function ${\bf
\phi}_{\omega}: \pi_2(M)\rightarrow\R$, $A\mapsto\int_A\omega$,
induced a {\it weight homomorphism} ${\bf \phi}:\Gamma\rightarrow\R$
which is injective. Then Hofer-Salamon \cite{HS} showed that the
group $\Gamma$ is isomorphic to a free abelian group with finite
many generators. We suppose $\{e_1,\cdots,e_k\}$ be the basis of
$\Gamma$, so for any $A\in\Gamma$, we have $A=\sum_{i=1}^k A_ie_i$.
Let $t=(t_1,\cdots,t_k)$, then we denote $t^A$ for $\prod_{i=1}^k
t_i^{A_i}$. Then the Novikov ring $\Lambda_{\omega,\phi}$ is the set
of formal sums as $$\lambda=\sum_{A\in\Gamma}\lambda_At^A,$$ where
$\lambda_A\in\Q$ such that for any constant $c>0$
$$\#\{ A\in\Gamma|\ \lambda_A\neq 0,\ {\bf \phi}(A)\le c \}<\infty.$$
Since the coefficient is the rational field $\Q$, our Novikov ring
is also a field. We note that the multiplication in Novikov ring is
$$\lambda\cdot\mu=\sum_{A,B\in\Gamma}\lambda_A\cdot\mu_Bt^{A+B}.$$

It is easy to verify that the following defined {\it scalar product}
$$\lambda\cdot\xi=\sum_{\tilde{x}\in Crit(F)}(\sum_{A\in\Gamma}\lambda_A\cdot
\xi_{(-A)\#\tilde{x}})\tilde{x}$$ is still in $C_*$, where the
connect sum $(-A)\#\tilde{x}$ is induced from the $\pi_2(M)$ action
on $\tilde{\Omega}_{\phi}$. Then we can consider the space $C_*$ as
a finite dimensional vector space over the field
$\Lambda_{\omega,\phi}$ with dimension of $\#{\rm Fix}(\phi_H)$.

Since from the first part of Corollary \ref{COR} we know that when
$\mu_{rel}(\tilde{y},\tilde{x})=\mu_{CZ}(\tilde{x})-\mu_{CZ}(\tilde{y})=1$,
the number $\#(P{\Cal M}^\nu(\tilde{y},\tilde{x}))$ is finite, then
we just define the boundary operator $\delta_\nu: C_*\rightarrow
C_*$ as: for any $\tilde{x}\in Crit_n(F)$
$$\delta_\nu(\tilde{x})=\sum_{\tilde{y}\in Crit_{n-1}(F)}\#(P{\Cal M}^\nu(\tilde{y},\tilde{x}))\tilde{y}.$$
Then by the conclusion in the second part of Corollary \ref{COR}, we
know that for any two $\tilde{z},\tilde{x}$ with
$\mu_{rel}(\tilde{z},\tilde{x})=2$, the oriented number
$$\#[\partial(P{\Cal
M}^\nu(\tilde{z},\tilde{x}))]=\sum_{\mu_{rel}(\tilde{z},\tilde{y})=1;\mu_{rel}(\tilde{y},\tilde{x})=1}
\#(P{\Cal M}^\nu(\tilde{z},\tilde{y}))\times\#(P{\Cal
M}^\nu(\tilde{y},\tilde{x}))=0.$$ So we have for any $\tilde{x}\in
Crit_n(F)$,
$$\delta_\nu^2(\tilde{x})=\sum_{\mu_{CZ}(\tilde{z})=n-2}[\sum_{\mu_{CZ}(\tilde{y})=n-1}
\#(P{\Cal M}^\nu(\tilde{z},\tilde{y}))\times\#(P{\Cal
M}^\nu(\tilde{y},\tilde{x}))]\tilde{z}=0.$$

Thus, we just define the Floer homology associated to $(J,H,\nu)$ of
the symplectic manifold $(M,\omega)$ and a symplectomorphism $\phi$
as the homology of the chain complex $(C_*,\delta_\nu)$, denoted by
$FH_*(J,H,\phi,\nu)$ or $FH_*(J,H,\phi)$.

As the last step, we should prove that the Floer homology groups
$FH_*(J,H,\phi,\nu)$ are independent of the almost complex structure
$J_t$ and the time-dependent Hamiltonian $H$ used to define them. To
do this, we need a continuation argument as the standard method used
in \cite{F3}\cite{SZ}\cite{LT}, {\it etc.}.

Given a fixed symplectomorphism $\phi$ and two suitable triples
$(J_0,H_0,\nu_0)$ and $(J_1,H_1,\nu_1)$ which are fine in the
construction as above. We want to show that
$$FH_*(J_0,H_0,\nu_0)\cong FH_*(J_1,H_1,\nu_1).$$ Thus, we need to
show there exists a chain homotopy. Firstly, we define a chain
homomorphism $$\Phi^0_1:
(C_*(H_0),\delta_{J_0,H_0,\nu_0})\rightarrow(C_*(H_1),\delta_{J_1,H_1,\nu_1}).$$

Suppose we have a family of generic pairs $(J_s,H_s)$, $s\in\R$, so
that $(J_s,H_s)\equiv(J_0,H_0)$ for $s\le 0$ and
$(J_s,H_s)\equiv(J_1,H_1)$ for $s\geq 1$. For two critical points in
different spaces $$\tilde{x}_0\in \widetilde{{\rm
Fix}}(\phi_{H_0})\subset\widetilde{\Omega}_\phi(H_0)\ \ {\rm and}\ \
\tilde{x}_1\in \widetilde{{\rm
Fix}}(\phi_{H_1})\subset\widetilde{\Omega}_\phi(H_1),$$ we can
similarly as before define the moduli space $$P{\Cal
M}(J_s,H_s,\tilde{x}_0,\tilde{x}_1)$$ of stable continuation
trajectories, which consists of stable $(J_s,H_s)$-orbits connecting
$\tilde{x}_0$ and $\tilde{x}_1$, say the element is
$$V=((v_1,\cdots,v_K),(f_1,\cdots,f_l),o): \Sigma\rightarrow M,\
v_j=u_j|_{\R\times[0,1]},$$ with some differences in that on each
main component $\Sigma_m$, the map $u_m:\ \Sigma_m\rightarrow M$,
$m=1,\cdots,K$, satisfies the following equation
\begin{equation}\label{010}
\bar{\partial}_{J_s,H_s}(u_m)=\frac{\partial u_m}{\partial s}(s,t)+
J_{s,t}(u_m(s,t))(\frac{\partial u_m}{\partial
t}(s,t)-X_{t}(H_s))=0.
\end{equation}
And we require that there is a $m_0\in\{1,\cdots,K\}$ so that when
$m<m_0$, $u_m$ is a stable $(J_0,H_0)$-map and when $m>m_0$ it is a
stable $(J_1,H_1)$-map satisfying the equation (\ref{001}).

We can apply all above construction to this new stable moduli space,
for example, we can similarly define the ambient space ${\Cal
B(J_s,H_s,\tilde{x}_0,\tilde{x}_1)}$, the neighborhood $\Cal W$ and
the compact virtual moduli space $P{\Cal
M}^{\nu_s}(J_s,H_s,\tilde{x}_0,\tilde{x}_1)$, {\it etc.}. Since for
those stable $(J_s,H_s)$-orbits the $s$-invariance does not hold for
the distinct main component $u_{m_0}$, the dimension of $P{\Cal
M}^{\nu_s}(\tilde{x}_0,\tilde{x}_1)$ will be
$$\mu_{rel}(\tilde{x}_0,\tilde{x}_1)=\mu_{CZ}(H_1,\tilde{x}_1)-\mu_{CZ}(H_0,\tilde{x}_0).$$

So when $\mu_{CZ}(H_1,\tilde{x}_1)=\mu_{CZ}(H_0,\tilde{x}_0)$,
$P{\Cal M}^{\nu_s}(\tilde{x}_0,\tilde{x}_1)$ is a finite set with
well-defined rational oriented number $\#(P{\Cal
M}^{\nu_s}(\tilde{x}_0,\tilde{x}_1))$. Consequently, we can define
the homomorphism between the spaces with same grade
$C_n(H_0)\rightarrow C_n(H_1)$ (then it is also a homomorphism
between $C_*(H_0)$ and $C_*(H_1))$ as
$$\Phi^0_1(\tilde{x}_0)=\sum_{\mu_{CZ}(\tilde{x}_1)=n}\#(P{\Cal
M}^{\nu_s}(J_s,H_s,\tilde{x}_0,\tilde{x}_1))\tilde{x}_1,$$ for
$\tilde{x}_0\in C_n(H_0)$.

Similarly, we can define a chain homomorphism $$\Phi^1_0:
(C_*(H_1),\delta_{J_1,H_1,\nu_1})\rightarrow(C_*(H_0),\delta_{J_0,H_0,\nu_0}).$$

Then just as the classical method of introducing an extra parameter
$\rho$ and applying all constructions as before to the two
parameters family $(J_s^\rho,H_s^\rho)$, used by Floer
\cite{F1}-\cite{F3} with the modification in that we replace the
classical moduli space of $(J_s,H_s)$-orbits ${\Cal M}(J_s,H_s)$ by
the branched labeled pseudomanifold, {\it i.e} the virtual moduli
space $P{\Cal M}^{\nu_s}(J_s^\rho,H_s^\rho)$ and correspondingly we
replace ${\Cal M}(J_s^\rho,H_s^\rho)$ by $P{\Cal
M}^{\nu_s^\rho}(J_s^\rho,H_s^\rho)$, we can prove that there exist
chain homotopies $$\Phi^1_0\circ\Phi^0_1\sim Id_{C_*(H_0)}\ \ {\rm
and} \ \ \Phi^0_1\circ\Phi^1_0\sim Id_{C_*(H_1)}.$$ Therefore, we
have induced an isomorphism $$(\Phi^0_1)^*:
FH_*(J_0,H_0,\nu_0)\rightarrow FH_*(J_1,H_1,\nu_1).$$ We omit the
details, and refer the reader to
\cite{F3}\cite{LT}\cite{HS}\cite{SZ},{\it etc.} for similar
arguments with respect to establishing isomorphism of Floer
(co)homologies for Hamiltonian periodic solutions. So we can just
denote the Floer homology associated with a compact symplectic
manifold $(M,\omega)$ and a symplectomorphism $\phi$ by
$FH_*(M,\omega,\phi,\nu)$ or simply by $FH_*(\phi,\nu)$.

\bigskip

Although $FH_*(\phi,\nu)$ are in general dependent of the
symplectomorphism $\phi$, we will show that they only depend on
$\phi$ up to Hamiltonian isotopy. That is to say, there is a natural
isomorphism
$$FH_*(\phi_0)\rightarrow FH_*(\phi_1)$$ if $\phi_0$ and
$\phi_1$ are related by a Hamiltonian isotopy. Let
$$\phi_t=\varphi_t^{-1}\circ\phi_0$$ be a Hamiltonian isotopy  from $\phi_0 $ to
$\phi_1=\varphi_1^{-1}\circ\phi_0$. That is to say, there exist
$t$-dependent Hamiltonian function $G_t$ and Hamiltonian vector
field $Y_t$ on $M$ satisfying
$$\frac{d}{dt}\varphi_t=Y_t\circ\varphi_t,\ \ \ \iota(Y_t)\omega=dG_t.$$
Recall for symplectomorphism $\phi_0=\phi$ satisfying
$H_t=H_{t+1}\circ\phi_0$, we have a map $u:\R^2\rightarrow M$ which
is the solution of the equation (\ref{001}) with boundary condition
(\ref{002}) and limits (\ref{003}). Then we define a new map
$v:\R^2\rightarrow M$ so that
\begin{equation}
v(s,t)=\varphi_t^{-1}(u(s,t)),\ \ J'_t=\varphi_t^*J_t.
\end{equation}
We see that $\phi_1^*J'_{t+1}=J_t$ and $v(s,t)$ is the solution of
the following equation
\begin{equation}\label{012}
\frac{\partial v}{\partial s}+J'_t(v)(\frac{\partial v}{\partial
t}-\varphi_t^{-1}\circ(X_t-Y_t)\circ\varphi_t(v))=0
\end{equation}
with boundary condition
\begin{equation}\label{013}
v(s,t+1)=\phi_1(v(s,t)).
\end{equation}
and the solution has limits
\begin{equation}\label{014}
\lim_{s\rightarrow\pm\infty}v(s,t)=\varphi_t^{-1}\psi_t(x_{\pm}),\ \
x_{\pm}=\phi_{H}x_{\pm}),
\end{equation}
where recall $x_{\pm}$ are nondegenerate fixed points of
$\phi_H=\psi_1^{-1}\circ\phi$. Thus there is a one-to-one
correspondence between the solutions of (\ref{001}), (\ref{002}) and
the solutions of (\ref{012}), (\ref{013}). We can extend the
arguments above to the virtual case in a similar way to get a 1-1
correspondence between their related virtual moduli space $P{\Cal
M}^{\nu_0}$ and $P{\Cal M}^{\nu_1}$. So we can obtain the
isomorphism $FH_*(\phi_0)\simeq FH_*(\phi_1)$.

In particular, as what Dostoglou-Salamon claimed in \cite{DS}, by
the ``fixed point index" property of the Maslov index listed in
section 2, the Euler characteristic is just the Lefschetz number of
$\phi$
$$\chi(FH_*(\phi))=\sum_{x\in{\rm Fix(\phi_H)}}{\rm sign\
det}(Id-d\phi_H(x))=L(\phi).$$

\bigskip
\noindent {\it Remark.} Recall the construction of Floer homology is
related to the choices of those based paths $\gamma_0$ in each
connected component of $\Omega_{\phi}$. In fact, we should denote
the Floer homology by $FH_*(\phi,\gamma_0)$. Nevertheless, we can
show that our construction above depends only on the homotopy class
of those based paths. A similar argument can be found in \cite{FO3}.

\bigskip
\noindent\small Institute of Mathematical Science\\
\noindent\small Nanjing University\\
\noindent\small Nanjing 210093, P.R.China\\
\noindent\small E-mail:  hailongher@ims.nju.edu.cn

\end{document}